\documentclass{amsart}
\usepackage{amsmath}
\usepackage{amsfonts}
\usepackage{amssymb}
\usepackage{amscd}
\usepackage{verbatim,multicol}
\usepackage[arrow,matrix,cmtip]{xy}

\usepackage{fullpage}

\newcommand{\trdeg}{\operatorname{tr.deg}}

\newcommand{\End}{\operatorname{End}}

\newcommand{\reg}{\operatorname{reg}}

\newcommand{\spec}{\operatorname{Spec}}

\newcommand{\mc}{\mathcal}

\newcommand{\mb}{\mathbb}

\newcommand{\dra}{\dashrightarrow}

\newcommand{\GK}{\operatorname{GK}}

\newcommand{\wt}{\widetilde}

\newcommand{\aut}{\operatorname{Aut}}

\DeclareMathOperator{\Pic}{Pic} \DeclareMathOperator{\HB}{H}

\newcommand{\se}{$\sigma$-positive}
\DeclareMathOperator{\calHom}{\mathcal{H}\mathit{om}}

\numberwithin{equation}{section}
 \theoremstyle{plain}
\newtheorem{theorem}[equation]{Theorem}

\newtheorem{lemma}[equation]{Lemma}

\newtheorem{proposition}[equation]{Proposition}

\theoremstyle{definition}
\newtheorem{definition}[equation]{Definition}
\newtheorem{notation}[equation]{Notation}

\newtheorem{remark}[equation]{Remark}

\newtheorem{hypothesis}[equation]{Hypothesis}
\newtheorem{example}[equation]{Example}

\begin{document}

\title{GK-dimension of birationally commutative surfaces}

\author{D. Rogalski}
\address{Department of Mathematics, UCSD, La Jolla, CA 92093-0112, USA.}
\email{drogalsk@math.ucsd.edu}
\thanks{The author was  partially supported by the NSF  through
grants DMS-0202479  and  DMS-0600834.}

\keywords{GK-dimension, graded rings, noncommutative projective
geometry, noncommutative surfaces, birational geometry}
\subjclass[2000]{14A22, 14E05, 16P90, 16S38, 16W50}

\begin{abstract}  Let $k$ be an algebraically closed field, let
$K/k$ be a finitely generated field extension of transcendence
degree $2$ with automorphism $\sigma \in \aut_k(K)$, and let $A
\subseteq Q = K[t; \sigma]$ be an $\mb{N}$-graded subalgebra with
$\dim_k A_n < \infty$ for all $n \geq 0$.  Then if $A$ is big enough
in $Q$ in an appropriate sense, we prove that $\GK A = 3,4,5,$ or
$\infty$, with the exact value depending only on the geometric
properties of $\sigma$.  The proof uses techniques in the birational
geometry of surfaces which are of independent interest.
\end{abstract}

\maketitle
\tableofcontents

\section{Introduction}
Throughout this paper, let $k$ be an algebraically closed field. Let
$A = \bigoplus_{n \geq 0} A_n$ be a finitely generated
$\mb{N}$-graded $k$-algebra which is  \emph{locally finite} ($\dim_k
A_n < \infty$ for all $n \geq 0$).  There is a close relationship
between such graded rings which are commutative and projective
algebraic geometry; for example, a projective variety of dimension
$n$ arises as $\operatorname{Proj} A$ for a commutative graded
domain $A$ of Krull dimension $n+1$.  Speaking broadly, the subject
of noncommutative projective geometry attempts to find interesting
ways to generalize to noncommutative rings this correspondence
between commutative graded $k$-algebras and projective $k$-schemes.
The Gelfand-Kirillov (GK)-dimension is typically the most useful
notion of dimension for rings in this theory, though it does not
always match the geometric intuition. In any case, noncommutative
$\mb{N}$-graded domains of GK-dimension $2$, which correspond to
noncommutative projective curves, have been classified in geometric
terms by Artin and Stafford in \cite{AS}. The classification of
noncommutative projective surfaces is an important ongoing research
problem, but it is less clear precisely which rings should
correspond to surfaces. For our purposes, given a locally finite
$\mb{N}$-graded domain $A$ which has a graded quotient ring $Q(A)
\cong D[t, t^{-1}; \sigma]$, we want to consider $A$ as
corresponding to a surface if the division ring $D$ has
transcendence degree $2$ over $k$.  An important special case occurs
when $D = K$ is a field, in which case we call the domain $A$
\emph{birationally commutative}.  The goal of this paper is to
calculate the GK-dimension of birationally commutative surfaces.

The main tools to be used in the calculation are the theory of
commutative projective surfaces and the Riemann-Roch theorem. There
is a major difficulty to overcome:  automorphisms of fields of
transcendence degree $2$ may not correspond nicely to automorphisms
of projective surfaces.  More specifically, given a finitely
generated field extension $K/k$ and an automorphism $\sigma \in
\aut_k(K)$, we say that $\sigma$ is \emph{geometric} if there exists
a projective variety $X$ with $k(X) = K$ and an automorphism $\tau:
X \to X$ which induces by pullback of rational functions the
automorphism $\sigma$. Although if $\trdeg K/k = 1$ every
automorphism of $K$ is geometric, this sometimes fails for fields of
higher transcendence degree.  The main thrust of this article is to
develop techniques to deal with non-geometric automorphisms.

We are now ready to state our main result.  To finesse the minor
technicality that graded quotient rings do not always exist, we
consider more generally the GK-dimension of all locally finite
$\mb{N}$-graded subalgebras $A \subseteq K[t, t^{-1}; \sigma]$ which
are \emph{big} in an appropriate sense (see
Definition~\ref{bigdef}).

\begin{theorem}
\label{intro-thm} Let $k$ be algebraically closed, let $K/k$ be a
finitely generated field extension with $\trdeg K/k = 2$, and let
$\sigma \in \aut_k(K)$. Then every big locally finite
$\mb{N}$-graded subalgebra $A$ of $Q = K[t, t^{-1}; \sigma]$ has the
same GK-dimension $d \in \{ 3,4,5, \infty \}$. Moreover, if $d <
\infty$, then $d = 4$ if and only if $\sigma$ is not geometric, and
in case $d = \infty$ then $A$ has exponential growth.
\end{theorem}
In fact, the value of $\GK A$ in Theorem~\ref{intro-thm} is
determined by geometric data associated to the automorphism
$\sigma$; see Theorem~\ref{final-thm} below for a more specific
statement.

Given a projective $k$-scheme $X$, an automorphism $\sigma: X \to X$
and an invertible sheaf $\mc{L}$ on $X$, the \emph{twisted
homogeneous coordinate ring} $B = B(X, \mc{L}, \sigma)$ is the ring
$\bigoplus_{n \geq 0} \HB^0(X, \mc{L}_n)$, where $\mc{L}_n$ is
defined by $\mc{L} \otimes \sigma^*(\mc{L}) \otimes \dots \otimes
(\sigma^{n-1})^* \mc{L}$, and for $x \in B_m, y \in B_n$ the
multiplication is defined by $x \star y = x \otimes
(\sigma^m)^*(y)$.  Such rings have good properties (such as the
noetherian property) when $\mc{L}$ satisfies an additional condition
called $\sigma$-ampleness.   For geometric automorphisms $\sigma \in
\aut_k(K)$, Theorem~\ref{intro-thm} is already known by work of
Zhang and the author \cite[Theorem 1.6]{RZ}.  The method is as
follows: choose $X$ projective with $k(X) = K$ and with an
automorphism $\sigma: X \to X$ corresponding to $\sigma \in
\aut(K)$. Then for any big subalgebra $A \subseteq Q = K[t, t^{-1};
\sigma]$ one may show that $A$ is both contained in and contains a
twisted homogeneous coordinate ring $B$. In this way one obtains
$\GK A = \GK B(X, \mc{L}, \sigma)$ for an appropriate choice of
($\sigma$-ample) $\mc{L}$. Moreover, Artin and Van den Bergh showed
that for a surface $X$  and  $\sigma$-ample sheaf $\mc{L}$, the
value of $\GK B(X, \mc{L}, \sigma)$  is either $3, 5$ or $\infty$
\cite[Theorem 1.7]{AV}. In fact, for higher-dimensional $X$ Keeler
gave bounds on the possible values of $\GK B(X, \mc{L}, \sigma)$
\cite[Theorem 6.1]{Ke1}, but less is known about what values are
actually obtained.

For a non-geometric $\sigma \in \aut_k K$, we cannot pick a nice
model $X$ with corresponding automorphism as in the last paragraph,
so we make do by carefully picking some projective nonsingular
surface $X$ with $k(X) = K$, and working with the \emph{birational}
map $\sigma: X \dra X$ corresponding to $\sigma \in \aut(K)$.
Fortunately, there is recent work by Diller and Favre \cite{DF} on
the dynamics of birational self-maps of smooth surfaces over
$\mb{C}$, which contains many helpful ideas. If $N^1(X) = \Pic
X/\sim$ is the group of divisors on $X$ modulo numerical
equivalence, then there is a natural way to define an action
$\sigma^*: \Pic(X) \to \Pic(X)$ which is given roughly by pullback
of divisors, and which descends to an action $\sigma^*: N^1(X) \to
N^1(X)$.  In fact $N^1(X) \cong \mb{Z}^m$ for some $m$ and
$\sigma^*$ is given by some (possibly non-invertible) integer matrix
$P$.  As is the case in \cite{AV} and \cite{Ke1}, this matrix is the
key to our GK-dimension calculations.   However, the notion of
pullback by a birational map is not always well-behaved under
iteration; importantly, one can choose the nonsingular model $X$ of
$K$ so that the birational map $\sigma: X \dra X$ is also
\emph{stable}, in the sense that the pullback map $(\sigma^n)^*$ is
given by the matrix $P^n$ for all $n \geq 1$ \cite[Theorem 0.1]{DF}.
The cases appearing in Theorem~\ref{intro-thm} then correspond to
the possible Jordan forms of the matrix $P$, which are classified in
\cite{DF} (see Theorem~\ref{DF-sum-thm} below.)

Since the paper \cite{DF} is so fundamental to our study here, we
give a rather thorough review in \S\ref{DF-sec} below of those
results from that paper that we need.  Since the authors of
\cite{DF} work exclusively over $\mb{C}$ and use analytic methods
and terminology, one of our aims here is to notice that many of
their proofs use only standard algebraic surface theory, which is
valid for any algebraically closed field.  In a few cases where a
proof in \cite{DF} is highly  analytic, we offer an alternative
algebraic proof.  One of the useful aspects of the \cite{DF} theory
is that it shows that non-geometric automorphisms $\sigma: K \to K$
of fields $K$ of transcendence degree $2$ are rather constrained:
for some model $X$ of $K$, the corresponding birational map $\sigma:
X \dra X$ either preserves a rational fibration of $X$, or else the
induced action $P: N^1(X) \to N^1(X)$ has an eigenvalue of modulus
bigger than $1$.  We offer some observations of our own concerning
the property of geometricity in \S\ref{geom-sec}.

Given a surface $X$ with a stable birational map $\sigma: X \dra X$
and invertible sheaf $\mc{L}$, we can set $\mc{L}_n = \mc{L} \otimes
\sigma^*\mc{L} \otimes \dots \otimes (\sigma^{n-1})^* \mc{L}$ and
build an analog of a twisted homogeneous coordinate ring, which we
write as $\wt{B}(X, \mc{L}, \sigma) = \bigoplus_{n \geq 0} \HB^0(X,
\mc{L}_n)$ and call a \emph{twisted section ring}. This construction
reduces to the usual twisted homogeneous coordinate ring when
$\sigma$ is an automorphism.   We give the construction of these
rings in \S\ref{B-sec}. Unfortunately, for non-geometric $\sigma$
such rings are badly behaved:  if $\mc{L}$ is also ample enough in
some sense, then we show that $\wt{B}(X, \mc{L}, \sigma)$ is not
finitely generated as an algebra.  However, we can still calculate
the growth of the graded pieces of $\wt{B}$ using the Riemann-Roch
theorem, and we do this calculation in \S\ref{growth-sec}.  Given a
finitely generated big subalgebra $A$ of $K[t, t^{-1}; \sigma]$, it
is easy to show that $A$ is contained in some $\wt{B}$, which gives
an upper bound on $\GK A$.  Finding a lower bound on $\GK A$ is
trickier, though, since in general we cannot find a copy of some
$\wt{B}$ as a subalgebra of $A$.  For the lower bound, we restrict
all sections to some curve $E$ on $X$ and reduce the problem to some
growth estimates for divisors on the curve $E$. The details are
given in \S\ref{lower-sec}. Then the upper and lower bound estimates
are combined to prove the main theorem in \S\ref{final-sec}.

The main theorem has applications to the classification theory of
noncommutative projective surfaces. In \cite[Theorem 1.1]{RS1}, we
completely described $\mb{N}$-graded domains $A$ which are
noetherian, generated in degree $1$, and have $Q(A) = K[t,t^{-1};
\sigma]$ for a finitely generated field extension $K$ with $\trdeg
K/k = 2$ and $\sigma$ geometric. Namely, in large degree such an $A$
is isomorphic either to a twisted homogeneous coordinate ring $B(X,
\mc{L}, \sigma)$, or to a \emph{na{\"\i}ve blowup algebra} $R(X, Z,
\mc{L}, \sigma)$ as studied in \cite{KRS} and \cite{RS2}. For a
non-geometric $\sigma$, what finitely generated domains $A$ with
$Q(A) = K[t,t^{-1}; \sigma]$ might look like is still unclear; we
don't know, for example, if any such $A$ are noetherian.  Indeed, we
have developed here the background material about non-geometric
automorphisms and twisted section rings in greater detail than we
really need for the proof of the main theorem about GK-dimension,
because we hope this will be of use in further work to better
understand such domains $A$. In any case, Theorem~\ref{intro-thm}
gives us a first understanding of what the noetherian examples (if
any) for non-geometric $\sigma$ look like: they have GK-dimension
$4$ (since noetherian $\mb{N}$-graded algebras cannot have
exponential growth). Additionally, Theorem~\ref{intro-thm} allows
the hypothesis that $\sigma$ is geometric in \cite[Theorem 1.1]{RS1}
to be replaced with the alternative hypothesis that $\GK A =3$ or
$5$.

\section{Review of results of Diller and Favre}
\label{DF-sec}

In this section, the word \emph{surface} will always mean a
nonsingular integral projective surface over the algebraically
closed field $k$.  The letters $W, X, Y, Z$  will always stand for
surfaces.  On such a surface $X$, we can identify
$\operatorname{Pic} X$ with the set of Weil divisors $D$ up to
linear equivalence or with invertible sheaves up to isomorphism.
There is a symmetric nondegenerate bilinear intersection form $( \,
\cdot \, , \, \cdot\, )$ on $\operatorname{Pic} X$, as defined in
\cite[Section V.1]{Ha}. We say that  divisors $D, E$ are
\emph{numerically equivalent} if $(D-E.C) = 0$ for all divisors $C$.
The quotient of $\Pic X$ by the relation of numerical equivalence is
the \emph{Neron-Severi group} $N^1(X) = \Pic X/\sim$, which is a
free abelian group of finite rank. A divisor $D$ is \emph{nef} if
$(D.C) \geq 0$ for all irreducible curves $C$ on $X$. We also work
with $\mb{R}$-divisors up to numerical equivalence, that is, the
group $N^1(X)_{\mb{R}} = N^1(X) \otimes_{\mb{Z}} \mb{R}$; the
intersection form uniquely extends to a form on $N^1(X)_{\mb{R}}$,
and  $D \in N^1(X)_{\mb{R}}$ is called \emph{nef} if $(D.C) \geq 0$
for all irreducible curves $C$ on $X$. The set of nef classes forms
the \emph{nef cone} $\operatorname{Nef}(X) \subseteq
N^1(X)_{\mb{R}}$. The \emph{pseudoeffective cone}
$\overline{\operatorname{NE}}(X) \subseteq N^1(X)_{\mb{R}}$ is the
smallest closed cone in $N^1(X)_{\mb{R}}$ containing the classes of
all effective curves. The nef cone and pseudoeffective cone are dual
with respect to the intersection form, and a nef class is also
pseudoeffective. See \cite{La} for more information on all of these
standard facts.

We now review the notion of pullback of divisors by a birational
map.
\begin{definition}
\label{alt-def-rem} \label{pullback-def} Let $f: X \dra Y$ be a
birational map of  surfaces. The map $f$ is defined except at a
finite set of \emph{fundamental points}.  Let $S \subseteq X$ be any
finite set of points in $X$ containing all of the fundamental points
of $f$, and let $U = X \smallsetminus S$, so $f \vert_U: U \to Y$ is
a morphism. Then given any divisor $D$ on $Y$, the divisor $( f
\vert_U)^*(D)$ on $U$ extends uniquely to a divisor $f^*(D)$ on $X$.
In this way, we define a map $f^*: \Pic Y \to \Pic X$.
\end{definition}

In the proofs below, it is useful to have a more explicit definition
of the pullback map $f^*$ which reduces to the case of a monoidal
transformation, and we discuss this next.  Suppose first that $\rho:
W \to X$ is a birational surjective morphism. Then $\rho = \rho_1
\rho_2 \dots \rho_n$ is a composition of finitely many monoidal
transformations $\rho_i: W_i \to W_{i-1}$ for some surfaces $W_0 =
X, \, W_1,\, \dots, W_n = W$ \cite[Corollary V.5.4]{Ha}. If $E
\subseteq W_i$ is the exceptional curve for $\rho_i$, then set $F_i
= (\rho_{i+1} \rho_{i+2} \dots \rho_n)^* (E)$.  It then follows from
\cite[Propositions V.3.2, V.3.6]{Ha} that the pullback map $\rho^*:
\Pic X \to \Pic W$ is an injection, with

\begin{equation}
\label{Pic-eq} \Pic W = \rho^*(\Pic X) \oplus \mb{Z} F_1 \oplus
\dots \oplus \mb{Z} F_n,
\end{equation}
where $(F_i.F_i) = -1$ for all $i$, $(F_i.F_j) = 0$ for $i \neq j$,
$(D.F_i) = 0$ for all $D \in \rho^*(\Pic X)$ and all $i$, and
$(\rho^* C_1, \rho^* C_2)_W = (C_1.C_2)_X$ for all $C_1, C_2 \in
\Pic X$. Moreover, each $F_i$ is a nonnegative linear combination of
the finitely many irreducible curves contracting under $\rho$ and so
is effective. Given any irreducible curve $C \subseteq X$, if
$\wt{C} \subseteq W$ is the proper transform of $C$ then iterating
\cite[Proposition V.3.6]{Ha} we have $\rho^*(C) = \wt{C} + \sum a_i
F_i$ for some nonnegative $a_i$.  We can also pushforward divisors:
$\rho_*: \Pic W \to \Pic X$ simply projects onto the factor
$\rho^*(\Pic X) \cong \Pic X$. Alternatively, it is clear that given
an irreducible curve $D \subseteq Y$, then $\rho_*(D) = \rho(D)$ if
$\rho(D)$ is a curve or $0$ if $\rho(D)$ is a point, and extending
this rule linearly determines the map $\rho_*$.

Now if $f: X \dra Y$ is a birational map, then there exists another
surface $W$ and birational morphisms $\phi: W \to X$, $\psi: W \to
Y$ such that $f = \psi \phi^{-1}$ \cite[Theorem V.5.5]{Ha}. Then it
is easy to check that $f^* = \phi_* \psi^*: \Pic Y \to \Pic X$ is
the same map as defined in Definition~\ref{pullback-def}; in
particular, this is independent of the choice of $W$.  We also
define pushforward of divisors by $f$ by setting $f_* = (f^{-1})^*
=\psi_* \phi^*: \Pic X \to \Pic Y$. The pullback map $f^*$ respects
numerical equivalence, and so it also induces maps $f^*: N^1(Y) \to
N^1(X)$ and $f^*: N^1(Y)_{\mb{R}} \to N^1(X)_{\mb{R}}$.

In the paper \cite{DF}, Diller and Favre define the pullback of
cohomology classes by a birational map $f: X \dra Y$ in case $k =
\mb{C}$, using the language of complex geometry \cite[Definition
1.8]{DF}.  They then give a very useful and beautiful classification
of birational self-maps of surfaces.  Our main aim in this section
is to review their theory and to show that a large number of their
results hold over any algebraically closed field. For most of their
proofs, we simply checked that upon substituting the algebraic
language of divisors for their more analytic language, the same
proof works.   In this case, we simply state their result below
without proof.   In a few cases where the translation of their proof
to the algebraic setting is non-obvious, or where we need a slightly
stronger formulation of a result, we provide a proof here. We review
not only the main results from \cite{DF} we need in the sequel, but
also the subsidiary results on which they depend. We hope this will
help the careful reader who wants to check that everything goes
through over arbitrary algebraically closed fields, as well as give
all readers a better idea of the flavor of the theory.

We begin with several results which follow very formally from the
definitions of pullback and pushforward.

\begin{lemma}\cite[Proposition 1.11(2)(3)]{DF}
\label{adjoint-lem} Let $f: X \dra Y$ be a birational map. Then
\begin{enumerate}
\item (Adjointness)  $(f^*C.D)_X = (C.f_*D)_Y$ for all $C \in
N^1(Y)_{\mb{R}}$, $D \in N^1(X)_{\mb{R}}$.

\item The map $f^* : N^1(Y)_{\mb{R}} \to N^1(X)_{\mb{R}}$
preserves the nef and pseudoeffective cones, in other words
$f^*(\operatorname{Nef}(Y)) \subseteq \operatorname{Nef}(X)$ and
$f^*(\overline{\operatorname{NE}}(Y)) \subseteq
\overline{\operatorname{NE}}(X)$.
\end{enumerate}
\end{lemma}

Pullback by a birational map does not generally respect the
intersection form. The next result calculates the exact nature of
the discrepancy.
\begin{lemma}\cite[Corollary 3.4]{DF}
\label{push-pull-lem} Let $f: X \dra Y$ be a birational map and
choose another surface $W$ and birational morphisms $\phi: W \to X$,
$\psi: W \to Y$ such that $f = \psi \phi^{-1}$. Decompose $\Pic W
\cong \phi^*(\Pic X) \oplus \mb{Z} F_1 \oplus \dots \oplus \mb{Z}
F_n$ with respect to the morphism $\phi$, as in \eqref{Pic-eq}, and
let $E_i = \psi_*(F_i) \in \Pic Y$. Then each $E_i$ is a sum with
nonnegative coefficients of irreducible curves contracting under
$f^{-1}$.  Moreover:
\begin{enumerate}

\item $(f^*(C).f^*(D)) = (C.D) + \sum_{i = 1}^n (C.E_i)(D.E_i)$ for
all $C, D \in N^1(Y)_{\mb{R}}$.

\item $(f^*(D).f^*(D)) = (D.D)$ for some $D \in N^1(Y)_{\mb{R}}$ if and only if $(D.V_i) = 0$ for all
irreducible curves $V_i \subseteq Y$ contracted by $f^{-1}$.
\end{enumerate}
\end{lemma}
\begin{proof}
The proof of \cite[Corollary 3.4]{DF}, which needs no essential
change, depends only on the earlier results \cite[Proposition 3.1,
Proposition 3.2, Theorem 3.3]{DF}.  These results also go through
with only obvious changes.  In fact, the algebraic versions of
\cite[Proposition 3.1 and Proposition 3.2]{DF} are immediate from
\cite[Proposition V.3.6 and Corollary V.3.7]{Ha}  (note that we
assume $f$ is birational whereas $f$ is allowed to be any surjective
morphism in these results in \cite{DF}; in particular, the constant
$\lambda_2$ of \cite[Proposition 3.1]{DF} satisfies $\lambda_2 = 1$
in our setting.)
\end{proof}

The behavior of the pullback maps under composition also needs
careful examination.

\begin{lemma}\cite[Proposition 1.13]{DF}
\label{comp-lem} Suppose that $f: X \dra Y$ and $g: Y \dra Z$ are
birational maps, and consider $f^*: \Pic Y \to \Pic X$, $g^*: \Pic Z
\to \Pic Y$, and $(gf)^*: \Pic Z \to \Pic X$.
\begin{enumerate}
\item  For any nef divisor $D \in \Pic Z$, $f^{*} g^{*}(D) -
(gf)^{*}(D)$ is effective.

\item If there does not exist any curve $C \subseteq X$ such that
$f(C) = p$ is a fundamental point of $g$, then $f^{*} g^{*} =
(gf)^{*}$.
\end{enumerate}
\end{lemma}
\begin{proof}
Since the proof in \cite{DF} seems quite analytic, we provide an
algebraic proof here. We may find a surface $W_1$ and birational
morphisms $\phi: W_1 \to X$, $\psi: W_1 \to Y$ such that $f = \psi
\phi^{-1}$, and a surface $W_2$ and birational morphisms $\mu: W_2
\to Y$, $\nu: W_2 \to Z$ such that $g = \nu \mu^{-1}$. Then let $h =
\mu^{-1} \psi: W_1 \dra W_2$ and find a surface $W_3$ and birational
morphisms $\tau: W_3 \to W_1$, $\theta: W_3 \to W_2$ such that $h =
\theta \tau^{-1}$.

Let $D \in \Pic Z$, and let $E = \nu^* D \in \Pic W_2$.  Let $P$ be
the finite set of fundamental points of $\mu^{-1}$, and let $Q = \{x
\in W_1 | \psi(x) \in P \}$, which is a proper closed subset of
$W_1$. Let $G_1, G_2, \dots, G_m \subseteq W_1$ be the distinct
irreducible curves contained in $Q$. Let $T = W_1 \setminus Q$ and
$U = Y \setminus P$. Then the birational map $\mu^{-1} \psi$ is
represented by the morphism $(\mu^{-1} \vert_U) \circ (\psi
\vert_T): T \to W_2$. Similarly, let $R \subseteq W_1$ be the finite
set of fundamental points of $\tau^{-1}$ and put $V = W_1 \setminus
R$, so the birational map $\theta \tau^{-1}$ is represented by the
morphism $\theta \circ (\tau^{-1} \vert_V): V \to W_2$.  Since
$\mu^{-1} \psi = \theta \tau^{-1}$ as birational maps, it now
follows from Definition~\ref{alt-def-rem} that $\psi^* \mu_* E -
\tau_* \theta^* E$ must be supported along the curves $G_1, \dots ,
G_m$, i.e. $\psi^* \mu_* E - \tau_* \theta^* E = \sum_{i =1}^m a_i
G_i \in \Pic W_1$.  Then $f^* g^* = \phi_* \psi^* \mu_* \nu^*$ and
$(gf)^* = \phi_* \tau_* \theta^* \nu^*$ and so $f^* g^*(D) -
(gf)^*(D) = \sum_{i = 1}^m a_i \phi_*(G_i)$.  In particular, in the
case of part (2), the set of $G_i$'s is empty and so the equation
$f^{*} g^{*} = (gf)^{*}$ follows.

Now suppose that $D$ is nef as in part (1).  Then $E = \nu^* D$ is
also nef by Lemma~\ref{adjoint-lem}(2). Write $\Pic W_1 =
\psi^*(\Pic Y) \bigoplus \mb{Z} F_1 \bigoplus \dots \bigoplus \mb{Z}
F_n$ with respect to the morphism $\psi$ as in \eqref{Pic-eq}. Since
$\psi^* \mu_* E \in \psi^*(\Pic Y)$, $(\psi^* \mu_* E.F_i) = 0$ for
all $i$. Because $\tau_* \theta^* E$ is also nef by
Lemma~\ref{adjoint-lem}(2) and $F_i$ is effective, $(\tau_* \theta^*
E.F_i) \geq 0$ for each $i$.  Setting $B = \psi^* \mu_* E - \tau_*
\theta^* E$, we have $(B.F_i) \leq 0$ for all $i$. Since each $G_i$
contracts under $\psi$, $G_i \in \sum \mb{Z} F_i$. Also, by the
first part of the proof, we have $B \in \sum \mb{Z} G_i$. Thus $B =
\sum_{i= 1}^n b_i F_i$.  Since $(F_i.F_j) = - \delta_{ij}$, this
forces $b_i \geq 0$ for all $i$ and so $B$ is effective. Finally,
$f^{*} g^{*}(D) - (gf)^{*}(D) = \phi_*(B)$ is then also effective.
\end{proof}

Now we want to concentrate on the special case of a birational map
$\sigma: X \dra X$ from a surface $X$ to itself.  In particular, we
want to study the pullback map $\sigma^*: N^1(X)_{\mb{R}} \to
N^1(X)_{\mb{R}}$. Recalling that $N^1(X)$ is a free Abelian group
$\cong \mb{Z}^d$ for some $d$, the map $\sigma^*$ is given by some
matrix in $M_d(\mb{Z})$.
\begin{example}
\label{cremona-ex} As Lemma~\ref{comp-lem} suggests, the pullback
map $\sigma^*$ may behave unpredictably with respect to iteration.
For example, suppose that $X = \mb{P}^2$ and $\sigma: \mb{P}^2 \dra
\mb{P}^2$ is the Cremona map defined by $(a:b:c) \mapsto
(bc:ac:ba)$. Then $\Pic X = N^1(X) \cong \mb{Z}$ and $\sigma^{*}:
\mb{Z} \to \mb{Z}$ is multiplication by $2$.  However, $\sigma^2$ is
the identity and so $(\sigma^2)^{*}$ is also the identity.
\end{example}

The previous example shows that we cannot expect $(\sigma^n)^{*} =
(\sigma^{*})^n$ to hold in general.  The next definition singles out
birational maps without this deficiency.
\begin{definition}
The birational map $\sigma: X \dra X$ is called \emph{stable} (the
term \emph{analytically stable} is used in \cite{DF}) if there does
not exist an irreducible curve $C \subseteq X$ and $n \geq 1$ such
that $\sigma^n(C)$ is a fundamental point of $\sigma$.
\end{definition}
For example, the Cremona map $\sigma: \mb{P}^2 \dra \mb{P}^2$
described above is not stable, because $\sigma$ contracts the line
$\{(a:b:c) | a = 0 \}$ to the point $(1:0:0)$ where $\sigma$ is
undefined.
\begin{remark}
\label{stable-rem} It is an easy consequence of the definition,
using that a point $p \in Y$ is a fundamental point for some
birational map $\rho: Y \dra Z$ if and only if $\rho^{-1}$ contracts
some curve to $p$  \cite[p. 256]{Sh}, that $\sigma: X \dra X$ is
stable if and only if $\sigma^{-1}$ is stable.  Also, note that if
$\sigma: X \dra X$ is stable and $C$ is a curve such that $\sigma(C)
= p_0$ is a point, then $\sigma(p_0) = p_1$ is defined; then since
$\sigma^2(C) = p_1$, $\sigma(p_1) = p_2$ is defined, and so on.
Continuing inductively, we can define a sequence of points $p_i$ by
the rule $\sigma(p_i) = p_{i+1}$.
\end{remark}

\begin{lemma}\cite[Theorem 1.14]{DF}
\label{stable-prop-lem} If $\sigma: X \dra X$ is stable then
$(\sigma^n)^{*} = (\sigma^{*})^n$ and $(\sigma^n)_* = (\sigma_*)^n$
for all $n \geq 1$.
\end{lemma}
\begin{proof}
This is an immediate consequence of Lemma~\ref{comp-lem}(2) and
Remark~\ref{stable-rem}.
\end{proof}

If $\sigma: X \dra X$ is a birational map and $f: Y \dra X$ is
another birational map, then we say that the birational map $\tau =
f^{-1} \sigma f : Y \dra Y$ is \emph{conjugate} to $\sigma$.  One of
the most fundamental results in \cite{DF} is that any birational
self-map is conjugate to a stable one.
\begin{theorem}\cite[Theorem 0.1]{DF}
\label{stable-choice-thm} Let $\sigma: X \dra X$ be a birational
map.  Then there is another surface $Y$ and a birational morphism
$\pi:Y \to X$ such that $ \pi^{-1} \sigma \pi: Y \dra Y$ is stable.
\end{theorem}

Next, for a birational map $\sigma: X \dra X$ we study how the maps
$(\sigma^n)^*: N^1(X)_{\mb{R}} \to N^1(X)_{\mb{R}}$ grow with $n$.
In this paper, it is most convenient to measure the growth of
functions using a different equivalence relation than that
traditionally used in the theory of GK-dimension \cite[p.5]{KL}.
Given two functions $f, g : \mb{N} \to \mb{R}$ which are monotone
increasing and positive valued for all $n \gg 0$, we write $f
\preccurlyeq g$ if there is a constant $d
> 0$ such that $f(n) \leq d g(n)$ for all $n \gg 0$.  If $f
\preccurlyeq g$ and $g \preccurlyeq f$ then we write $f \sim g$.
There is no difference between this definition and the one given in
\cite[p.5]{KL} for functions of polynomial growth, but our
definition distinguishes between various functions of exponential
growth such as $f(n) = n 2^n$ and $g(n) = 2^n$.  It is also useful
to note that if $\lim_{n \to \infty} f(n)/g(n)$ is a finite positive
number, then $f \sim g$.

For any surface $X$, $|| \cdot ||$ will indicate some arbitrary
matrix norm on the space
$\operatorname{End}_{\mb{C}}(N^1(X)_{\mb{C}})$. The particular
choice of matrix norm will never impact any of our results below.
Recall that $\rho = \max \{ | \lambda | \big| \lambda \in \mb{C}\
\text{is an eigenvalue of}\ \sigma^* \}$ is called the
\emph{spectral radius} of $\sigma^*$.

\begin{lemma} \cite[Lemma 1.12]{DF}
\label{rho-facts-lem} Let $\sigma: X \dra X$ be birational map and
let $\sigma^*: N^1(X)_{\mb{R}} \to N^1(X)_{\mb{R}}$ have spectral
radius $\rho$. Then there is a nef class $D \in N^1(X)_{\mb{R}}$
which is an eigenvector for $\sigma^*$ with eigenvalue $\rho$.
\end{lemma}
\begin{proof}
We note here that this result follows quickly from existing
references.  By Lemma~\ref{adjoint-lem}, $\sigma^*$ preserves the
nef cone $\operatorname{Nef}(X) \subseteq N^1(X)_{\mb{R}}$. Then the
spectral radius $\rho$ is an eigenvalue of $\sigma^*$ with an
eigenvector in $\operatorname{Nef}(X)$ by \cite[Theorem 3.1]{V}.
(This is the same argument used by Keeler, for the case where
$\sigma$ is an automorphism, in \cite[Lemma 3.2]{Ke1}; only the fact
that $\sigma^*$ preserves the nef cone is needed.)
\end{proof}

\begin{lemma}
\label{growth-inv-lem} \label{power-growth-lem} \cite[Corollary
1.16]{DF} Let $\sigma: X \dra X$ be a birational map.
\begin{enumerate}
\item For any birational map $\tau: Y \dra Y$ which is conjugate to
$\sigma$,  $|| (\sigma^n)^* || \sim  || (\tau^n)^* ||$.
\item $|| (\sigma^n)^* || \sim n^j \rho^n$ for some uniquely determined real
$\rho \geq 1$ and integer $j \geq 0$.
\item If $\sigma$ is stable,
then $\rho$ is the spectral radius of $\sigma^*$ and $j+1$ is the
size of the largest Jordan block in $\sigma^*$ associated to an
eigenvalue of modulus $\rho$.
\end{enumerate}
\end{lemma}
\begin{proof}
The proof of \cite[Corollary 1.16]{DF}, which depends on
\cite[Proposition 1.15]{DF} and Lemma~\ref{comp-lem}(1), goes
through with little change to prove part (1).  If $\sigma$ is
stable, then $(\sigma^n)^* = (\sigma^*)^n$.  Choosing a basis of
$N^1(X)_{\mb{C}}$ so that $\sigma^*$ is in Jordan form, each
coordinate of the matrix $(\sigma^*)^n$ is a $\mb{C}$-linear
combination of terms of the form $n^a \lambda^i$ where $a \geq 0$,
$0 \leq i \leq n$, and $\lambda$ is an eigenvalue of $\sigma^*$.  If
$j+1$ is the size of the largest Jordan block in $\sigma^*$
associated to an eigenvalue $\lambda$ of modulus $\rho$, then $n^j
\lambda^n$ is the fastest growing such term; thus $|| (\sigma^n)^*
|| \sim n^j \rho^n$ in this case, proving part (3). Since each
$(\sigma^n)^*$ is in $\operatorname{End}_{\mb{Z}}(N^1(X)) \cong
M_d(\mb{Z})$ and there is a positive lower bound on the norms of
nonzero integer matrices, it follows that $\rho \geq 1$.  For a
$\sigma$ which is not necessarily stable, it now follows from
Theorem~\ref{stable-choice-thm} and part (1) that $|| (\sigma^n)^*
|| \sim n^j \rho^n$ for some $\rho,j$ with $\rho \geq 1$.   The
values of $\rho, j$ are obviously uniquely determined given our
definition of equivalence of functions, completing the proof of part
(2).
\end{proof}

\begin{definition}
\label{bir-data-def}  Given a birational map $\sigma: X \dra X$,  we
associate to $\sigma$ the \emph{growth data} $(\rho, j)$ as
determined by Lemma~\ref{growth-inv-lem}(2).
\end{definition}

The main achievement of \cite{DF} is to classify the possible values
of growth data $(\rho, j)$ associated to birational maps $\sigma: X
\dra X$ and to describe which kinds of surfaces can appear in each
case. The last results from \cite{DF} which we recall in this
section contain the main work for this classification.  We need a
few more definitions. First, a \emph{rational fibration} of $X$ is a
morphism $f: X \to C$, where $C$ is a nonsingular projective curve
and the generic fiber of $f$ is isomorphic to $\mb{P}^1$. A
birational map $\sigma: X \dra X$ \emph{preserves} such a fibration
if there is an automorphism $\tau: C \to C$ such that $ f \sigma =
\tau f$ (as birational maps).

\begin{lemma}\cite[Proposition 1.7]{DF}
\label{nonpos-int-lem} Let $\sigma: X \dra X$ be a birational map
which is not an automorphism, with associated growth data $(\rho,
j)$.  We choose a surface $W$ and morphisms $\phi: W \to X$, $\psi:
W \to X$ such that $\sigma = \psi \phi^{-1}$; moreover we make this
choice such that the rank of $N^1(W)$ is minimal. Choose any
decomposition $\phi = \phi_1 \phi_2 \dots \phi_m$, where each
$\phi_i$ is a single monoidal transformation, and let $E \subseteq
W$ be the exceptional curve contracted by $\phi_m$. Then $V =
\psi(E) \subseteq X$ is an irreducible curve which is contracted by
$\sigma^{-1}$, and $(V.V) \geq -1$.  Moreover, we have the
following.
\begin{enumerate}
\item  If $(V.V) = -1$, then there exists a monoidal transformation $\pi:
X \to X'$ with exceptional curve $V$ (and in this case, if $\sigma$
is stable then so is $\sigma' = \pi \sigma \pi^{-1}: X' \dra X'$.)

\item If $(V.V) = 0$, then there exists a
rational fibration $f: X \to C$ such that $V \cong \mb{P}^1$ is a
generic fiber.
\end{enumerate}
\end{lemma}
\begin{proof}
The same proof as that of  \cite[Proposition 1.7]{DF} works here; we
have just organized the result differently by moving some parts of
the proof into the statement.    The proof of the existence of a
rational fibration in case (2) depends on the classical result
\cite[Proposition 4.3]{BPV}, which also is valid over an arbitrary
algebraically closed field.  Note that \cite[Proposition 1.7]{DF}
contains an extra hypothesis that $(C.C) \leq 0$ for all irreducible
curves $C$ contracting under $\sigma^{-1}$ (such as $V$).  This
hypothesis is not needed in our version; we simply draw no
conclusion if $(V.V) > 0$.
\end{proof}

\begin{proposition}
\label{small-growth-prop}
 \cite[Lemma 4.1]{DF}
Suppose that $\sigma: X \dra X$ is a birational map, such that $\|
(\sigma^n)^* \|$ is bounded.  Then $\sigma$ is conjugate to an
automorphism $\tau: Y \to Y$ such that $(\tau^n)^*: N^1(Y) \to
N^1(Y)$ is the identity for some $n \geq 1$.
\end{proposition}
\begin{proof}
The proof in \cite{DF}, which is easy to rewrite in algebraic
language, uses Theorem~\ref{stable-choice-thm},
Lemma~\ref{growth-inv-lem}, Lemma~\ref{push-pull-lem}(2), and
Lemma~\ref{nonpos-int-lem}.
\end{proof}

\begin{proposition} \cite[Lemma 4.2]{DF}
\label{big-growth-prop} Let $\sigma: X \dra X$ be birational, with
growth data $(\rho, j)$ given by Definition~\ref{bir-data-def} such
that $\rho = 1$.  Then either (1) $\sigma$ is conjugate to an
automorphism $\tau: Y \to Y$, or else (2) $j = 1$ and $\sigma$ is
conjugate to another birational map $\mu: Z \dra Z$ such that (i)
$\mu$ is stable; (ii) $Z$ has a rational fibration $f: Z \to C$
which is preserved by $\mu$;  and (iii) if $S$ is the set of
irreducible curves contracted by $\mu^{-1}$, then for all $V \in S$
we have $\mu^*(V) = V$ in $N^1(Z)$, and $(V_i. V_{\ell}) = 0$ for
all $V_i, V_{\ell} \in S$.
\end{proposition}
\begin{proof}
This is a restatement and slight strengthening of the result of
\cite[Lemma 4.2]{DF}, which does not require in case (2) the
restrictions given in (iii).  So we indicate the needed adjustments
to the proof.

First, by Theorem~\ref{stable-choice-thm} we replace $X$ with a
blowup to assume $\sigma$ is stable.   Suppose that $\sigma$ is not
an automorphism.  As in the statement of Lemma~\ref{nonpos-int-lem},
choose a $W$ with Neron-Severi group of minimal rank which has
morphisms $\phi: W \to X$, $\psi: W \to X$ such that $\sigma = \psi
\phi^{-1}$.  Given any point $p \in X$ which is a fundamental point
for $\sigma$, we can choose a decomposition $\phi = \phi_1 \phi_2
\cdots \phi_m$ where each $\phi_i$ is a single monoidal
transformation, and moreover such that the exceptional curve $E
\subseteq W$ contracted by $\phi_m$ satisfies $\phi(E) = p$.  Let $V
= \psi(E)$.  If $(V.V) = -1$, then using
Lemma~\ref{nonpos-int-lem}(1) we contract $V$ by a single monoidal
transformation $\pi: X \to X'$ such that $\sigma' = \pi \sigma
\pi^{-1}: X' \dra X'$ is again stable. Replacing $X$ by $X'$ and
$\sigma$ by $\sigma'$, we repeat this process as many times as
possible.  The process is finite since the rank of the Neron-Severi
group decreases by one with each contraction.

If the process of the last paragraph results in an automorphism, we
are done, so assume not.  Then we have obtained a stable birational
map $\sigma: X \dra X$ such that the set $S$ of irreducible curves
contracted by $\sigma^{-1}$ is nonempty. We can pick some
fundamental point $p$ of $\sigma$ and repeat the first part of the
previous paragraph; then $(V.V) \geq 0$ by
Lemma~\ref{nonpos-int-lem}, since we repeated the contraction
process until no longer possible.  Since $\rho = 1$, the proof of
\cite[Lemma 4.2]{DF} (which depends on Lemma~\ref{rho-facts-lem},
Lemma~\ref{push-pull-lem}(2), and Lemma~\ref{nonpos-int-lem}) shows
that $(V.V) = 0$, there is a rational fibration $f: X \to C$
preserved by $\sigma$ with generic fiber $V$, that $\sigma^*(V) = V$
in $N^1(X)$, and $j = 1$.   The union of the curves which are
contracted by $\sigma^{-1}$ to the point $p$ is a connected set by
Zariski's main theorem, and all of these curves lie in the same
fiber of $f$ as $V$ since $\sigma$ preserves the fibration.  But $V$
is a generic fiber, and so $V$ is the only curve contracted by
$\sigma^{-1}$ to $p$. Moreover, any other $V' \in S$ must also lie
in some fiber and so $(V'.V) = 0$.

Applying the previous paragraph to each fundamental point of
$\sigma$, we see that each $V_i \in S$ contracts by $\sigma^{-1}$ to
a unique point $p_i$, that $\sigma^*(V_i) = V_i$ in $N^1(X)$, and
that  $(V_i.V_{\ell}) = 0$ for all $V_{\ell} \in S$.
\end{proof}

\section{Geometricity}
\label{geom-sec}

A surface always means a nonsingular integral projective surface
over the algebraically closed field $k$ in this section. In
addition, $\sigma: X \dra X$ continues to denote a birational map of
a surface $X$.
\begin{definition}
If $K$ is a finitely generated field extension of $k$ with $\trdeg
K/k = 2$ and $\tau: K \to K$ is a field automorphism over $k$, we
say that $\tau$ is \emph{geometric} if there exists a surface $Z$
with $k(Z) = K$ such that the corresponding birational map $\tau: Z
\dra Z$ is an automorphism. Given a birational map $\sigma: X \dra
X$, we say that $\sigma$ is \emph{geometric} if the corresponding
field automorphism $\sigma: K \to K$ is geometric, or equivalently
if $\sigma$ is conjugate (as defined in the previous section) to
some automorphism $\tau: Z \to Z$ of a surface $Z$.
\end{definition}

In \cite{DF}, $\sigma: X \dra X$ is called \emph{bimeromorphically
conjugate to an automorphism} when we call it geometric.  The
following is a summary theorem of the main classification result in
\cite{DF}.

\begin{theorem}
\label{DF-sum-thm} \cite[Theorem 0.2]{DF}  Let $\sigma: X \dra X$ be
a birational map, and define the growth data $(\rho, j)$ by
Definition~\ref{bir-data-def}.  Choose some matrix norm $|| \cdot
||$ on $\End_{\mb{C}}(N^1(X)_{\mb{C}})$.  Then exactly one of the
following cases occurs:
\begin{enumerate}
\item $\rho = 1$, $j = 0$, $|| (\sigma^n)^* || \sim 1$, and $\sigma$ is
geometric; in fact, $\sigma$ is conjugate to some automorphism
$\tau: Y \to Y$ such that $(\tau^n)^*: N^1(Y) \to N^1(Y)$ is the
identity for some $n \geq 1$.
\item $\rho = 1$, $j = 1$, $|| (\sigma^n)^* || \sim n$, $\sigma$ is not geometric,
and $\sigma$ is conjugate to a birational map $\tau: Y \dra Y$ where
$Y$ has a rational fibration $f: Y \to C$ preserved by $\tau$.
\item $\rho = 1$, $j = 2$, $|| (\sigma^n)^* || \sim n^2$, and $\sigma$ is geometric.
\item $\rho > 1$, and $|| (\sigma^n)^* || \sim n^j \rho^n$ is growing
exponentially.
\end{enumerate}
\end{theorem}
\begin{proof}
Most of this follows immediately from
Propositions~\ref{small-growth-prop} and \ref{big-growth-prop} and
Lemma~\ref{growth-inv-lem}. It only remains to show that if $\rho =
1$ and $\sigma$ is geometric, then $j = 0$ or $2$. But in this case
we pass to a conjugate automorphism $\tau: Y \to Y$, and then the
fact that $j = 0$ or $2$ is shown by Artin and Van den Bergh
\cite[Lemma 5.4]{AV}.  For $k = \mb{C}$, this is also shown in
\cite[Theorem 4.3, Proposition 4.4]{DF} (relying on previous work of
Bellon, Cantat, and Gizatullin).  In fact these results show that in
case (3), it is even true that $\sigma$ preserves an elliptic
fibration.  We certainly expect that this also holds over an
arbitrary algebraically closed field, though we have not checked
since we will not need in the sequel that $\sigma$ preserves an
elliptic fibration in case (3).
\end{proof}

Theorem~\ref{DF-sum-thm} completely determines whether a birational
map $\sigma: X \dra X$ with $\rho = 1$ is geometric, in terms of the
value of $j$; in particular, non-geometric such maps have a very
restricted form.  Proposition~\ref{geom-char-prop} below offers
another characterization of geometricity, independent of the value
of $\rho$.  This gives a nice intuitive picture and will also be
useful in our study of twisted section rings in the next section. In
the following definition, we single out a pathology that turns out
to characterize stable non-geometric maps.
\begin{definition}
\label{self-con-def} Let $\sigma: X \dra X$ be birational. Suppose
that there exists a point $p \in X$ such that $\sigma^{-1}$ is
defined at $p_n = \sigma^{-n}(p)$ for all $n \geq 0$, but $\sigma^n$
is not defined at $p$ for infinitely many $n > 0$. Then we say that
$p$ is an \emph{unbalanced point} and that $\sigma$ is
\emph{unbalanced}. If no such $p$ exists then we say $\sigma$ is
\emph{balanced}.
\end{definition}
\begin{remark}
\label{unbal-rem}  The following easy observations, which are left
to the reader, will be useful below.  Let $\sigma: X \dra X$ be
birational.  If $\sigma$ is unbalanced, then  $\sigma^n: X \dra X$
is also unbalanced for every $n \geq 1$.  Conversely, if $\sigma$ is
stable and $\sigma^n$ is unbalanced for some $n \geq 1$, then
$\sigma$ is unbalanced.
\end{remark}

\begin{proposition} \label{geom-char-prop}
Let $\sigma: X \dra X$ be a birational map.
\begin{enumerate}
\item If $\sigma$ is unbalanced, then $\sigma$ is not geometric.

\item If $\sigma$ is stable, then $\sigma$ is balanced if and only if $\sigma$ is geometric.
\end{enumerate}
\end{proposition}
\begin{proof}
(1) Let $\sigma: X \dra X$ be unbalanced.   As a first step, we will
prove the following claim:  suppose that $\phi: Z \to X$ is a
birational morphism, and let $\tau = \phi^{-1} \sigma \phi: Z \dra
Z$; then $\tau^n: Z \dra Z$ is also unbalanced for some $n \geq 1$.
Since $\phi$ is a composition of finitely many monoidal
transformations, by induction it is enough to prove the claim in
case $\phi$ is itself a single monoidal transformation, say with
exceptional curve $E \subseteq Z$ and contracted point $q \in X$.
Let $p \in X$ be an unbalanced point for $\sigma$, and define $p_n =
\sigma^{-n}(p)$ for $n \geq 0$.   If $p_n \neq q$ for all $n \geq
0$, then define $r = \phi^{-1}(p)$, while if $p_n = q$ for exactly
one $n$, then define $r = \phi^{-1}(p_{n+1})$.  In either case it is
now easy to check that $r$ is an unbalanced point for $\tau$, so
$\tau$ itself is unbalanced.   Otherwise, $p_i = p_j = q$ for some
$0 < i < j$.   Let $d = j -i$; then $\sigma^{-d}$ is defined at $q$
with $\sigma^{-d}(q) = q$.  Now suppose that $\sigma^{-d}$ is a
local isomorphism at $q$. Then the inverse $\sigma^{d}$ is also
defined at $q$, and since $\sigma^n(q)$ is defined for all $n \leq
0$, it follows that $\sigma^n(q)$ is defined for all $n \in \mb{Z}$.
But then $\sigma^n(p)$ is defined for all $n \in \mb{Z}$, a
contradiction. Thus $\sigma^{-d}$ is not a local isomorphism at $q$,
so $\sigma^{-d}$ induces a non-surjective map of tangent spaces $T_q
\to T_q$.  But because points of the exceptional curve $E$ of $\phi$
correspond to tangent directions at $q$, it follows that
$\tau^{-d}(E) = r$ for some point $r \in E$.  Moreover, since
$\sigma^{-d}$ is defined at $q$, $\tau^{-d}$ is defined at every
point of $E$ (for example, by the universal property of blowing up);
in particular, at $r$.  Thus $r$ is an unbalanced point for
$\tau^d$, proving the claim.

Next, we claim that it is impossible to have an unbalanced
birational map $\tau: Z \dra Z$ and a birational morphism $\psi: Z
\to Y$ with $\mu = \psi \tau \psi^{-1}: Y \to Y$ an automorphism.
Suppose this is possible, and choose such an example where $\psi$ is
a composition of the minimal possible number $m$ of monoidal
transformations (obviously $m > 0$). Now let $q \in Z$ be an
unbalanced point. Then $\tau^n$ fails to be defined at $q$ for
infinitely many $n > 0$, say for some sequence $0 < n_1 < n_2 <
\dots$.  Thus we can find irreducible curves $C_1, C_2, \dots
\subseteq Z$ such that $\tau^{-n_i}(C_i) = q$.  Clearly we must have
that $\psi(C_i) = r_i$ is a point for all $i$.  Since $\psi$
contracts finitely many curves, we have $C_i = C_j$ for some $0 < i
< j$, and so $r_i = r_j = r$, say.  Moreover, $\mu^{-n_i}(r) =
\psi(q) = \mu^{-n_2}(r)$; letting $d = n_2 - n_1$, we then have
$\mu^{d}(r) = r$.  Work now instead with $\tau' = \tau^d: Z \dra Z$
and $\mu' = \mu^d: Y \to Y$. Then $\mu'(r) = r$, so if $\pi: \wt{Y}
\to Y$ is the monoidal transformation blowing up $r$, then $\mu'$
lifts to an automorphism $\nu = \pi^{-1} \mu' \pi: \wt{Y} \to
\wt{Y}$. Now $\psi = \pi \theta$ for some morphism $\theta: Z \to
\wt{Y}$, by \cite[Proposition V.5.3]{Ha}.  But $\tau'$ is again
unbalanced (Remark~\ref{unbal-rem}), and descends to an automorphism
$\nu = \theta \tau' \theta^{-1}: \wt{Y} \to \wt{Y}$, where $\theta$
is composed of fewer than $m$ monoidal transformations.  This
contradicts the minimal choice of $m$ and this contradiction proves
the claim.

Finally, suppose that $\sigma: X \dra X$ is unbalanced and
geometric. Then we may find a surface $Z$ and birational morphisms
$\phi: Z \to X$, $\psi: Z \to Y$, such that $\mu = \psi \phi^{-1}
\sigma \phi \psi^{-1}: Y \to Y$ is an automorphism.  By the first
claim, $\tau^n = \phi^{-1} \sigma^n \phi: Z \to Z$ is unbalanced for
some $n \geq 1$. Then by the second claim, $\mu^n$ is not an
automorphism, a contradiction.

\medskip

(2).  Let $\sigma: X \dra X$ be balanced and stable. Suppose that
$\sigma$ is not already an automorphism (in which case we are done)
and let $C \subseteq X$ be an irreducible curve with $\sigma(C) = p$
a point.  By the stability condition, $\sigma^n(p)$ is defined for
all $n \geq 0$, so at a point $r \in C$ for which $\sigma(r)$ is
defined, we have $\sigma^n(r)$ is defined for all $n \geq 0$. There
are at most finitely many points of $C$ at which $\sigma$ is not
defined; for one of these, say $q$, there must be a curve $D$ with
$\sigma^{-1}(D) = q$; then by the stability condition, $q_n =
\sigma^{-n}(q)$ is defined for all $n \geq 0$, with $\sigma^{-1}$
defined at each $q_n$, implying that $\sigma^n(q)$ is defined for
all $n \gg 0$ by the balanced hypothesis.  In conclusion, $\sigma^n$
is defined at every point of $C$ for some $n \geq 1$. Then we can
choose a nonsingular surface $Z$ and morphisms $\phi: Z \to X, \psi:
Z \to X$ for which $\psi \phi^{-1} = \sigma^{n}$, and where
$\phi^{-1}$ is a local isomorphism at every point of $C$ (see the
proof of \cite[Theorem V.5.5]{Ha}.)  Thus $(C.C)_X =
(\phi^{-1}(C).\phi^{-1}(C))_Z \leq -1$ since $\phi^{-1}(C)$ is a
curve contracting under $\psi$.  In conclusion, all irreducible
curves of $X$ contracted by $\sigma$ have negative
self-intersection.  Now using Lemma~\ref{nonpos-int-lem} (applied to
$\sigma^{-1}$) we can find a monoidal transformation $\pi: X \to X'$
blowing down some irreducible curve contracted by $\sigma$, where
replacing $X$ by $X'$ and $\sigma$ by $\sigma' = \pi \sigma
\pi^{-1}$, $\sigma'$ is still stable.  Moreover, $\sigma'$ is also
balanced, by the first claim proved in the first half of the proof,
together with Remark~\ref{unbal-rem}.

Repeat the previous paragraph as many times as possible.  It must
reach an automorphism after finitely many steps since the rank of
the Neron-Severi group decreases by $1$ each time.  Thus we have
shown that for a stable map, balanced implies geometric, and the
converse was proved in part (1).
\end{proof}

The preceding proposition is useful for constructing examples of
non-geometric maps.
\begin{example}
\label{non-geom-ex} We thank Michael Artin for suggesting the
following example.  Assume that $k$ is uncountable.  Let $X =
\mb{P}^2$, let $\tau: \mb{P}^2 \dra \mb{P}^2$ be the Cremona map of
Example~\ref{cremona-ex} and let $\phi: \mb{P}^2 \to \mb{P}^2$ be an
automorphism. Let $\sigma = \tau \phi: \mb{P}^2 \dra \mb{P}^2$. Let
$\ell_1, \ell_2, \ell_3$ be the coordinate lines which are
contracted by $\tau$. Note that $\sigma^{-1}(\ell_j) = p_{j,0}$ is a
point for each $j$.  Now if we choose $\phi$ generically (outside
some countable union of proper closed subsets of
$\operatorname{PGL}(2,k)$), we can ensure that $\sigma^n(\ell_j)$ is
a curve for all $j$ and $n \geq 1$, and moreover that we can
inductively define $p_{j, i+1} = \sigma^{-1}(p_{j,i})$ for all $j$
and $i \geq 0$.  Then $\sigma^{-1}$ and thus $\sigma$ is stable, and
$\sigma$ is not geometric by Proposition~\ref{geom-char-prop} since
each point $p_{j,0}$ is unbalanced.   Note also that $\Pic X =
\mb{Z}$ and $\sigma^*: \Pic X \to \Pic X$ is multiplication by $2$.
So $\rho = 2 > 1$ for this stable birational map.
\end{example}

We remark that there also geometric birational maps $\sigma: X \dra
X$ with $\rho > 1$; an automorphism of a K3 surface with $\rho
> 1$ is given in \cite[Example 3.9]{Ke1}.  Also, it is easy to find examples
of case (2) of Theorem~\ref{DF-sum-thm}; in fact, any ruled surface
has a non-geometric birational map preserving the ruling which is of
this type (see \cite[Remark 7.3]{DF}).

As a further application of Proposition~\ref{geom-char-prop}, we
have the following result that the property of geometricity is
invariant under base extension.
\begin{lemma}
\label{base-change-lem} Let $k \subseteq L$ where $L$ is another
algebraically closed field, let $K/k$ be a finitely generated field
extension with $\trdeg K/k = 2$, and let $\sigma \in \aut_k K$. Let
$F$ be the field of fractions of $K \otimes_k L$, so $\sigma$
extends uniquely to an automorphism $\wt{\sigma} \in \aut_L F$. Then
$\sigma \in \aut_k K$ is geometric if and only if $\wt{\sigma} \in
\aut_L F$ is geometric.
\end{lemma}
\begin{proof}
It is a standard result that since $k$ is algebraically closed, $K
\otimes_k L$ is a domain, so it makes sense to consider its field of
fractions $F$.  If $\sigma \in \aut_k K$ is geometric, pick a
projective model $X$ of $K$ such that the induced map $\sigma: X \to
X$ is an automorphism.  Then $\sigma$ lifts to an automorphism
$\wt{\sigma}: \wt{X} \to \wt{X}$ of the projective $L$-surface
$\wt{X} = X \times_{\spec k} \spec L$, so the corresponding field
automorphism $\wt{\sigma} \in \aut_L F$ is again geometric.

Conversely, if $\sigma \in \aut_k K$ is not geometric, we choose a
nonsingular projective model $X$ of $K$ such that the corresponding
birational map $\sigma: X \dra X$ is stable, by
Theorem~\ref{stable-choice-thm}. By
Proposition~\ref{geom-char-prop}(2), there is some unbalanced point
$p \in X$ for $\sigma$; then it is easy to see that the lifted point
$\wt{p} \in \wt{X} = X \times_{\spec k} \spec L$ is an unbalanced
point for $\wt{\sigma}: \wt{X} \dra \wt{X}$, so by
Proposition~\ref{geom-char-prop}(1), $\wt{\sigma} \in \aut_L F$ is
not geometric either.
\end{proof}

\section{Twisted section rings}
\label{B-sec}

If $\sigma: X \to X$ is an automorphism of a surface $X$, then given
an invertible sheaf $\mc{L}$ on $X$ we may form the \emph{twisted
homogeneous coordinate ring} $B(X, \mc{L}, \sigma)$ as described in
the introduction.  Since $B_n = \HB^0(X, \mc{L}_n)$ is the full
vector space of global sections of an invertible sheaf on $X$, many
questions about $B$ reduce to questions about invertible sheaves.
For example, the calculation of the growth of $\dim_k B_n$ can make
use of the Riemann-Roch theorem.  In this section, we generalize
this construction to form a ring $\wt{B}(X, \mc{L}, \sigma)$ where
$\sigma: X \dra X$ is just a stable birational map, and where
$\wt{B}_n = \HB^0(X, \mc{L}_n)$ is the global sections of the
analogous invertible sheaf $\mc{L}_n$ on $X$. The growth of the
graded pieces of $\wt{B}$ will still be a consequence of the
Riemann-Roch theorem, and we give it in the next section.
Unfortunately, for non-geometric $\sigma$ the ring $\wt{B}$ is
otherwise not well-behaved: we will show in this case that (for
$\mc{L}$ ample enough) $\wt{B}$ is not finitely generated as an
algebra.

We now describe some notational conventions that will hold
throughout this section.  As in previous sections, $X$ will always
be a nonsingular integral projective surface.  Given an invertible
sheaf $\mc{L}$ on $X$, for convenience we will assume that we have a
fixed embedding $\mc{L} \subseteq \mc{K}$ of $\mc{L}$ in the
constant sheaf $\mc{K}$ of rational functions $K = k(X)$; there is
then a corresponding Weil divisor $D$ on $X$ with $\mc{L} =
\mc{O}_X(D)$. Given some stable birational map $\sigma: X \dra X$,
we define the pullback map $\sigma^*: \Pic X \to \Pic X$ on Weil
divisors as in Definition~\ref{pullback-def}.   Thinking in terms of
invertible sheaves, we define $\sigma^* \mc{L} =
\mc{O}_X(\sigma^*(D))$, so there is also a fixed embedding $\sigma^*
\mc{L} \subseteq \mc{K}$.  As in the theory of twisted homogeneous
coordinate rings, it is convenient to use the notation
$\mc{L}^{\sigma} = \sigma^* \mc{L}$, and we define $\mc{L}_0 =
\mc{O}_X$ and $\mc{L}_n = \mc{L} \otimes \mc{L}^{\sigma} \otimes
\dots \otimes \mc{L}^{\sigma^{n-1}} \subseteq \mc{K}$ for $n \geq
1$. The birational map $\sigma: X \dra X$ induces an automorphism $K
\to K$ which we also call $\sigma$. Finally, we write $V^{\sigma} =
\sigma(V)$ for any subset $V \subseteq K$.

We begin by recording some simple sheaf-theoretic properties of the
pullback map.
\begin{lemma}
\label{pullback-sec-lem} Let $\sigma: X \dra X$ be a stable
birational map, and $\mc{L}$ an invertible sheaf on $X$.  Choose a
nonsingular surface $W$ with birational morphisms $\phi: W \to X$
and $\psi: W \to X$ such that $\sigma =  \psi \phi^{-1}$.
\begin{enumerate}
\item  $\mc{L}^{\sigma} \subseteq \mc{K}$ is the reflexive hull $\mc{F}^{**} = \calHom
(\calHom (\mc{F}, \mc{O}_X), \mc{O}_X)$ of the sheaf $\mc{F} =
\phi_* \psi^* \mc{L} \subseteq \mc{K}$.

\item  There is an induced injective pullback of sections map
\[
\sigma^*:  \HB^0(X, \mc{L}) \to \HB^0(X, \mc{L}^{\sigma})
\]
which is also the restriction of the automorphism $\sigma: K \to K$.

\item  If $V \subseteq \HB^0(X, \mc{L})$ generates $\mc{L}$, then
$V^{\sigma}$ generates $\mc{L}^{\sigma}$ except possibly at the
finitely many fundamental points of $\sigma$.  Similarly, $V
V^{\sigma} \dots V^{\sigma^{n-1}}$ generates $\mc{L}_n$ except at a
finite set of points for all $n \geq 1$.
\end{enumerate}
\end{lemma}
\begin{proof}
(1)  The sheaf $\mc{F} = \phi_* \psi^* \mc{L}$ is locally free
except possibly at the finite set $S$ of fundamental points of
$\sigma$. The reflexive hull $\mc{F}^{**}$ is the unique invertible
sheaf on $X$ agreeing with $\mc{F}$ on $U = X \smallsetminus S$
\cite[Sublemma 7.7]{RS1}.  Now part (1) is just a reinterpretation
of Definition~\ref{pullback-def}.

(2) Note that since we have embeddings $\mc{L} \subseteq \mc{K}$ and
$\mc{L}^{\sigma} \subseteq \mc{K}$, taking sections we also have
inclusions $\HB^0(X, \mc{L}) \subseteq K$ and $\HB^0(X,
\mc{L}^{\sigma}) \subseteq K$.   Now $\sigma^*: \HB^0(X, \mc{L}) \to
\HB^0(X, \mc{L}^{\sigma})$ may be defined by composing the injective
pullback map $\psi^*: \HB^0(X, \mc{L}) \to \HB^0(W, \psi^* \mc{L})$,
the bijective pushforward of sections $\phi_*: \HB^0(W, \psi^*
\mc{L}) \to \HB^0(X, \mc{F})$, and the injection $\HB^0(X, \mc{F})
\to \HB^0(X, \mc{F}^{**})$ induced by the canonical injection of
sheaves $\mc{F} \subseteq \mc{F}^{**} \subseteq \mc{K}$.  It follows
formally that $\sigma^*$ is also given by restricting $\sigma: K \to
K$.

(3) By the construction of $\sigma^*$ in part (2), it is clear that
$\sigma^*(V) = V^{\sigma} \subseteq \HB^0(X, \mc{F})$ generates the
sheaf $\mc{F}$ except possibly at the finite set $S$ of fundamental
points of $\sigma$. The injection $\mc{F} \subseteq \mc{F}^{**} =
\mc{L}^{\sigma}$ is an isomorphism at points not in $S$, so
$V^{\sigma}$ generates $\mc{L}^{\sigma}$ at points not in $S$.
 Applying this to each power $\sigma^i$, we see that $V^{\sigma^i}$
generates $\mc{L}^{\sigma^i}$ at points not in the set $S_i$ of
fundamental points of $\sigma^i$; thus $V V^{\sigma} \dots
V^{\sigma^{n-1}}$ generates $\mc{L}_n$ except possibly at points of
the finite set $\bigcup_{i = 1}^{n-1} S_i$.
\end{proof}

\begin{definition}
\label{wtb-def}
Given a stable birational map $\sigma: X \dra X$ and
$\mc{L} \subseteq \mc{K}$ an invertible sheaf, we define the
\emph{twisted section ring}
\[
\wt{B}(X, \mc{L}, \sigma) = \bigoplus_{n \geq 0} \HB^0(X,
\mc{L}_n)t^n \subseteq Q = K[t, t^{-1}; \sigma],
\]
with multiplication induced by that of $Q$.
\end{definition}

Note that if $\sigma$ is an automorphism, then $\wt{B}(X, \mc{L},
\sigma) \cong B(X, \mc{L}, \sigma)$ is the usual twisted homogeneous
coordinate ring.  Defining $\wt{B}$ as an explicit subring of $Q$ is
not really necessary, but it avoids a tedious proof that the
multiplication of $\wt{B}$ is associative.  Still, to ensure that
Definition~\ref{wtb-def} actually defines a ring, we need to check
that $\wt{B}$ is closed under multiplication.   For this, it is
enough to check that $\wt{B}_m \wt{B}_n \subseteq \wt{B}_{m+n}$ for
$m, n \geq 0$. If $f \in \HB^0(X, \mc{L}_m)$ and $g \in \HB^0(X,
\mc{L}_n)$, then $ft^m gt^n = f \sigma^m(g) t^{m+n}$ and
$\sigma^m(g) = (\sigma^m)^*(g) \in \HB^0(X, \mc{L}_n^{\sigma^m})$ by
part (2) of Lemma~\ref{pullback-sec-lem}. Now since $\sigma$ is
stable, we have
\begin{gather*}
\mc{L}_n^{\sigma^m} = (\sigma^m)^*(\mc{L} \otimes \sigma^* \mc{L}
\otimes \dots \otimes (\sigma^{n-1})^* \mc{L})  = (\sigma^m)^*
\mc{L} \otimes (\sigma^m)^* \sigma^*(\mc{L}) \otimes
\dots \otimes (\sigma^m)^* (\sigma^{n-1})^* \mc{L} \\
= (\sigma^m)^* \mc{L} \otimes (\sigma^{m+1})^*(\mc{L}) \otimes \dots
\otimes (\sigma^{n+m-1})^* \mc{L}
\end{gather*}
and so $\mc{L}_m \otimes \mc{L}_n^{\sigma^m} = \mc{L}_{m+n}$, as
subsheaves of $\mc{K}$.  Then $f \sigma^m(g)$ is the image of $f
\otimes (\sigma^m)^*(g)$ under the multiplication map $\HB^0(X,
\mc{L}_m) \otimes \HB^0(X, \mc{L}_n^{\sigma^m}) \to \HB^0(X,
\mc{L}_{m+n})$, so $\wt{B}$ is indeed a subalgebra of $Q$.

The goal of the rest of the section is to prove that the ring
$\wt{B}(X, \mc{L}, \sigma)$ is typically not finitely generated as a
$k$-algebra if $\sigma$ is non-geometric.  The exact theorem will
hold only for $\mc{L}$ ``ample enough".  We next define an
appropriate such class of invertible sheaves.  The definition is
fairly arbitrary; it is made just for the convenience of this paper
to include enough positivity properties for several later results.

\begin{definition}
\label{stab-eff-def} Given a stable birational map $\sigma: X \dra
X$, we say an  invertible sheaf $\mc{L}$ on $X$ is \emph{\se}\
 if (1) $\mc{L}$ is ample, (2) The sheaf
$\mc{L}_n$ is generated by its global sections for all $n \geq 1$,
and (3) $\HB^i(X, \mc{L}_n) = 0$ for all $i > 0$ and $n \geq 1$.
\end{definition}

We note that it is easy to find \se\ invertible sheaves.
\begin{lemma}
\label{stab-eff-lem} Given $\sigma: X \dra X$ stable and a very
ample invertible sheaf $\mc{M}$ on $X$, then $\mc{L} =
\mc{M}^{\otimes d}$ is \se\ for all $d \gg 0$.
\end{lemma}
\begin{proof}
The proof uses the notion of Castelnuovo-Mumford regularity with
respect to the very ample sheaf $\mc{M}$, and is mostly a matter of
quoting some known results. We recall the definition: for a coherent
sheaf $\mc{F}$ on $X$, $\mc{F}$ is called \emph{$d$-regular} if
$\HB^i(X, \mc{F} \otimes \mc{M}^{\otimes d-i}) = 0$ for all $i > 0$.
It is a fact that if $\mc{F}$ is $d$-regular, then $\mc{F}$ is
$e$-regular for all $e \geq d$ \cite[Theorem 1.8.5]{La}.  We let
$\reg \mc{F}$ be the smallest $d$ for which $\mc{F}$ is $d$-regular.

The formula $\lim_{n \to \infty} \reg \mc{M}^{\otimes n} = - \infty$
easily follows since $\mc{M}$ is ample.   By a theorem of Fujita,
there is a fixed bound $M$ such that $\reg \mc{F} \leq M$ for all
nef invertible sheaves $\mc{F}$ \cite[Theorem 1, p. 520]{Fj}.  Also,
there is a constant $C$ such that for any two invertible sheaves
$\mc{N}, \mc{P}$ on $X$, one has $\reg (\mc{N} \otimes \mc{P}) \leq
\reg \mc{N} + \reg \mc{P} + C$ \cite[Proposition 2.8]{Ke2}.

For $d \geq 1$ put $\mc{L} = \mc{M}^{\otimes d}$.  Then $\mc{L}$ is
nef and so $(\mc{L}^{\sigma} \otimes \dots \otimes
\mc{L}^{\sigma^{n-1}})$ is nef by Lemma~\ref{adjoint-lem}.  Then for
$d \gg 0$ we have $\reg \mc{L} \leq -C-M$ and so $\reg \mc{L}_n \leq
\reg \mc{L} + \reg (\mc{L}^{\sigma} \otimes \dots \otimes
\mc{L}^{\sigma^{n-1}})  + C \leq 0$.  Thus $\HB^i(X, \mc{L}_n) = 0$
for $i = 1, 2$.  Also, since $\mc{L}_n$ is $0$-regular, it is
generated by its global sections \cite[Theorem 1.8.5 and Remark
1.8.14]{La}.
\end{proof}

\begin{theorem}
\label{non-fg-thm} Let $\sigma: X \dra X$ be stable and
non-geometric, and suppose that $\mc{L} \subseteq \mc{K}$ is  \se.
Then $\wt{B}(X, \mc{L}, \sigma)$ is not a finitely generated
$k$-algebra.
\end{theorem}
\begin{proof}
Since $\sigma$ is not geometric, by Proposition~\ref{geom-char-prop}
there must exist an unbalanced point $q \in X$, so $q$ is a
fundamental point of $\sigma^i$ for infinitely many $i \geq 1$.  We
show first that in fact we can find a point which is a fundamental
point for all positive powers of $\sigma$.   Suppose we have $0 <
i_1 < i_2 < i_3$ such that $\sigma^{i_1}$ and $\sigma^{i_3}$ are not
defined at $q$, but $\sigma^{i_2}(q)$ is defined.  Then we can find
irreducible curves $F_1$ and $F_3$ such that $\sigma^{-i_1}(F_1) =
q$ and $\sigma^{-i_3}(F_3) = q$, so $\sigma^{i_2 -i_1}(F_1) =
\sigma^{i_2}(q) = \sigma^{i_2 -i_3}(F_3)$, and this contradicts the
fact that $\sigma$ is stable. Thus $\sigma^i$ is undefined at $q$
for all $i \gg 0$. Let $i_0$ be the largest nonnegative integer such
that $p = \sigma^{i_0}(q)$ is defined.  Then $p$ is a fundamental
point of $\sigma^i$ for all $i \geq 1$.  For each $i \geq 1$ we can
pick an irreducible curve $E_i \subseteq X$ such that
$\sigma^{-i}(E_i) = p$.

Now fix $n \geq 2$ and consider the multiplication map $\wt{B}_i
\otimes \wt{B}_{n-i} \to \wt{B}_n$ for some $0 < i < n$. Dropping
the powers of $t$, this may be identified with
\[
\theta_i: \HB^0(X, \mc{L}_i) \otimes \HB^0(X, \mc{L}_{n-i})
\overset{1 \otimes (\sigma^i)^*}{\longrightarrow}  \HB^0(X,
\mc{L}_i) \otimes \HB^0(X, \mc{L}_{n-i}^{\sigma^i}) \longrightarrow
\HB^0(X, \mc{L}_i \otimes \mc{L}_{n-i}^{\sigma^i})= \HB^0(X,
\mc{L}_n)
\]
where the second map is the natural multiplication map. We claim
that the image of $\theta_i$ is contained in $\HB^0(X, \mc{I}_p
\otimes \mc{L}_n)$ where $\mc{I}_p$ is the ideal sheaf of the point
$p$.  Suppose for the moment that claim has been proved for all $i$.
Then we will have that the image of $\bigoplus_{i = 1}^{n-1}
\wt{B}_i \otimes \wt{B}_{n-i} \to \wt{B}_n$ under the multiplication
map is contained in $\HB^0(X, \mc{I}_p \otimes \mc{L}_n) t^n$. Since
$\mc{L}$ is \se, $\mc{L}_n$ is generated by global sections at $p$
and so $\HB^0(X, \mc{I}_p \otimes \mc{L}_n) \subsetneq \HB^0(X,
\mc{L}_n)$.  It follows that $\wt{B}$ is not generated in degrees
$1$ through $n-1$, and since $n \geq 2$ was arbitrary the result
follows.

It remains to prove the claim.  For this, it will be enough to prove
that the image of the pullback of sections map $(\sigma^i)^*:
\HB^0(X, \mc{L}_{n-i}) \to \HB^0(X, \mc{L}_{n-i}^{\sigma^i})$ is
contained in $\HB^0(X, \mc{I}_p \otimes \mc{L}_{n-i}^{\sigma^i})$.
Let $\mc{N} =  (\mc{L}^{\sigma} \otimes \dots \otimes
\mc{L}^{\sigma^{n-i-1}})$ and note that $\mc{N}$ is nef by
Lemma~\ref{adjoint-lem}.   Then $\mc{L}_{n-i} = \mc{L} \otimes
\mc{N}$ is the product of an ample and a nef sheaf, so is also
ample.  Write $\mc{L}_{n-i} = \mc{O}_X(D)$ for some Weil divisor
$D$.  Given a section $0 \neq s \in \HB^0(X, \mc{L}_{n-i})$, let $D'
\in |D|$ be the divisor of zeroes of $s$, where $|D|$ is the
complete linear system of effective divisors linearly equivalent to
$D$. Note that since $\mc{L}_{n-i}$ is generated by its sections and
obviously $\mc{L}_{n-i} \not \cong \mc{O}_X$, we have $\dim_k
\HB^0(X, \mc{L}_{n-i}) \geq 2$.   It follows that for a generic $D'
\in |D|$, $D' = \sum  a_j C_j$ for some irreducible curves $C_j$
where (i) no $C_j$ contracts under $\sigma^{-i}$, and (ii) no $C_j$
contains a fundamental point of $\sigma^{-i}$.  We will prove that
for all such generic $D'$, the support of the divisor
$(\sigma^i)^*(D')$ contains $p$; then $(\sigma^i)^*(s) \in \HB^0(X,
\mc{I}_p \otimes \mc{L}_{n-i}^{\sigma^i})$ will hold for generic $s
\in \HB^0(X, \mc{L}_{n-i})$, so for all $s$ as required.  Thus, let
$D' \in |D|$ with $D' = \sum  a_j C_j$ satisfying conditions (i) and
(ii).  Since $D'$ is ample we have $(D'.E_i) > 0$. Furthermore, $C_j
\neq E_i$ for each $j$, so some $C_j$ intersects $E_i$ nontrivially;
but then $\sigma^{-i}$ is defined at every point of $C_j$, and
$\sigma^{-i}(C_j)$ is a curve in the support of $(\sigma^i)^*(D')$
which contains the point $\sigma^{-i}(E_i) = p$.
\end{proof}

The preceding theorem is not needed in the proof of the main theorem
of the paper, but we have included it as an interesting negative
result.  As part of the program to classify noncommutative
projective surfaces, one would like to describe in terms of geometry
all connected finitely generated $\mb{N}$-graded subalgebras $A$ of
$Q = K[t, t^{-1}; \sigma]$, where $K/k$ is a finitely generated
field extension of transcendence degree $2$.  For field
automorphisms $\sigma$ which are geometric, the paper \cite{RS1}
succeeds in this goal, for those $A$ which are also noetherian and
generated in degree $1$: such an $A$ is  equal (in large degree)
either to a twisted homogeneous coordinate ring $B(X, \mc{L},
\sigma)$ where $X$ is a surface, or else to a special kind of
subring $R(X, Z_{\mc{I}}, \mc{L}, \sigma) = \bigoplus_{n \geq 0}
\HB^0(X, \mc{I}_n \otimes \mc{L}_n)$ of such a $B$ called a
\emph{na{\"\i}ve blowup} as studied in \cite{KRS} and \cite{RS2}
(where here  $\mc{I}_n = \mc{I} \cdot \sigma^*\mc{I} \cdots
(\sigma^{n-1})^* \mc{I}$ for some ideal sheaf $\mc{I}$ defining a
$0$-dimensional subscheme $Z_{\mc{I}}$ of $X$). When $\sigma$ is
\emph{not} geometric, we do not know how to describe the subrings
$A$ of $Q$.  Theorem~\ref{non-fg-thm} shows that in this case the
natural analogs of twisted homogeneous coordinate rings, the
algebras $\wt{B}(X, \mc{L}, \sigma)$, unfortunately do not suffice
to describe any such $A$.   A different approach seems to be needed
to understand the non-geometric case; we hope to address this
question in future work.

\section{Growth of $\wt{B}$}
\label{growth-sec}

In this section, we calculate the growth of the graded pieces of a
ring $\wt{B}(X, \mc{L}, \sigma)$ using the Riemann-Roch formula. For
the case where $\sigma$ is non-geometric we will also need to rely
on the information from the classification result
Theorem~\ref{DF-sum-thm}.  The case where $\sigma$ is an
automorphism is already known (see \cite[Theorem 1.7]{AV} or
\cite[Theorem 6.1]{Ke1}), but we give a uniform proof that works in
all cases since this takes little extra effort.

The following is the situation we always consider from now on.
\begin{notation}
\label{standard-not} Let $\sigma \in \aut_k(K)$, where $K/k$ is a
finitely generated field extension of an algebraically closed field
$k$, with $\trdeg K/k = 2$. Choose a nonsingular integral projective
surface $X$ over $k$ with $k(X) = K$ such that the induced map
$\sigma: X \dra X$ is stable (which is possible by
Theorem~\ref{stable-choice-thm}.)  Fix some basis of $N^1(X) \cong
\mb{Z}^d$ and let $P \in M_d(\mb{Z})$ be the $d \times d$ matrix
representing the pullback map $\sigma^*: N^1(X) \to N^1(X)$. Fix any
matrix norm $\| \cdot \|$ on $M_d(\mb{C})$, so $|| P^n || \sim n^j
\rho^n$ where $(\rho, j)$ is the associated growth data (which is
independent of the choice of $X$ by Lemma~\ref{growth-inv-lem}). Let
$\mc{K}$ be the constant sheaf of rational functions on $X$.
\end{notation}

In some of our growth calculations, it will be convenient to assume
the following further restrictions on the choice of $X$:
\begin{hypothesis}
\label{standard-hyp}  Assume Notation~\ref{standard-not}.  If
$\sigma$ is geometric, then by definition we can  choose $X$ so that
the map $\sigma: X \to X$ is an automorphism. If $\sigma$ is
non-geometric with $\rho = 1$, then by
Proposition~\ref{big-growth-prop}(2) we have $j = 1$ and we can
choose $X$ so that $\sigma: X \dra X$ is stable; $\sigma$ preserves
some rational fibration $f: X \to C$; and if $S$ is the set of
irreducible curves contracted by $\sigma^{-1}$, then $P(V_i) = V_i$
in $N^1(X)$ and $(V_i.V_{\ell}) = 0$ for all $V_i, V_{\ell} \in S$.
\end{hypothesis}

We need the following simple linear algebra lemma.  The proof is
similar to the proof of Lemma~\ref{growth-inv-lem} and is left to
the reader.
\begin{lemma}
\label{Jordan-lem} Fix Notation~\ref{standard-not}, and let $Q \in
M_d(\mb{R})$ be an invertible matrix.
\begin{enumerate}
\item $\| Q P^n \| \sim n^j \rho^n$.

\item For any two column vectors $v, w \in \mb{R}^d$,  $ | v^{T} Q P^n w | \preccurlyeq
n^j \rho^n$.

\item There is a dense subset $U \subseteq \mb{R}^d$ (in fact, one can take $U$ to be the
complement of some quadric hypersurface) such that for any $v \in
U$, $| v^{T} Q P^n v | \sim  n^j \rho^n$.
\end{enumerate}
\end{lemma}

We now prove a series of growth estimates for the iterates of
divisor classes under the pullback map $P$.

\begin{lemma}
\label{form-growth-lem} Assume Notation~\ref{standard-not}.
\begin{enumerate}
\item For any $D, E \in N^1(X)_{\mb{R}}$, we have $| (P^nD.E) |
\preccurlyeq n^j \rho^n$.
\item For any ample $E \in N^1(X)_{\mb{R}}$, $ (P^nE.E)  \sim  n^j
\rho^n$.
\end{enumerate}
\end{lemma}
\begin{proof}
(1) For a review of the basic facts concerning intersection theory
on a surface, see \cite[Section V.1]{Ha}.  In particular, it is
standard that the intersection form on $X$ is given by some
symmetric real invertible matrix $Q$, such that for any column
vectors $D, E \in N^1(X)_{\mb{R}} \cong \mb{R}^d$ we have $(D.E) =
D^TQE$.  Then $(P^nD.E) = E^{T} Q P^n D$ and the result follows from
Lemma~\ref{Jordan-lem}(2).

(2) As in part (1), let $Q$ be the matrix of the intersection form.
Let $E$ be any ample divisor, and let $U$ be the open subset of
Lemma~\ref{Jordan-lem}(3).  Since the nef cone
$\operatorname{Nef}(X) \subseteq N^1(X)_{\mb{R}}$ is a cone of full
dimension $d$ in $\mb{R}^d$,  we can choose some $C \in U$ which is
nef.  We can also choose $m \gg 0$ so that $D = mE - C$ is nef. We
have
\begin{equation}
\label{form-growth-eq} (P^nC.C) = (P^n(mE-D).mE-D) = m^2(P^nE.E) -
(P^nD.C) - m(P^nE.D).
\end{equation}
Since $D$ and $C$ are nef by construction, and $P^nD, P^nE$, and
$P^nC$ are nef by Lemma~\ref{adjoint-lem}, \eqref{form-growth-eq}
implies that $m^2 (P^nE.E) \geq (P^nC.C) \geq 0$ for all $n$.  Using
that $C \in U$, we have $(P^nC.C) \sim n^j \rho^n$ by
Lemma~\ref{Jordan-lem}(3). Together with part (1) it follows that
$(P^nE.E) \sim n^j \rho^n$ as desired.
\end{proof}

\begin{lemma}
\label{form-growth-lem2} Assume that $\rho = 1$ and that
Hypothesis~\ref{standard-hyp} holds.  Let $D$ be an ample divisor on
$X$ and set $D_n = \sum_{i = 0}^{n-1} P^i(D)$. Then $(D_n.D_n) \sim
n^{j+2}$.
\end{lemma}
\begin{proof}

We claim first that for any $D \in N^1(X)$ there is an integer
constant $N \geq 0$ such that $(P^b(D). P^a(D)) = (P^{b-a}(D).D) +
aN$ for any $0 \leq a \leq b$.  If $\sigma$ is geometric, then we
have chosen $\sigma: X \to X$ be an automorphism in
Hypothesis~\ref{standard-hyp}.  In this case, $P$ is an invertible
matrix preserving the intersection form and so the claim obviously
follows, with $N = 0$.

Now suppose that $\sigma$ is non-geometric.  Then $j = 1$ and recall
that we have assumed that $X$ is chosen  in
Hypothesis~\ref{standard-hyp} so that if $S$ is the set of
irreducible curves contracted by $\sigma^{-1}$, then $P(V_i) = V_i$
in $N^1(X)$ and $(V_i.V_{\ell}) = 0$ for all $V_i, V_{\ell} \in S$.
Define the classes $E_i$ for the map $\sigma: X \dra X$ as in
Lemma~\ref{push-pull-lem}; each $E_i$ is a sum with nonnegative
coefficients of irreducible curves in $S$, and so $P(E_i) = E_i$ in
$N^1(X)$ and $(E_i.E_{\ell}) = 0$ for all $i,\ell$. Then by
Lemma~\ref{push-pull-lem}(1), we see that for any class $B \in
N^1(X)$, we have $(PB.E_i) = (PB.PE_i) = (B.E_i) + \sum_{\ell = 1}^d
(B.E_{\ell})(E_i.E_{\ell}) = (B.E_i)$.  Now for any $0 \leq a \leq
b$, by Lemma~\ref{push-pull-lem}(1), we have
\[
(P^bD.P^aD) = (P^{b-1}D.P^{a-1}D) + \sum_{i = 1}^d
(P^{b-1}D.E_i)(P^{a-1}D.E_i)
\]
and by the previous calculation we also have $\sum_{i = 1}^d
(P^{b-1}D.E_i)(P^{a-1}D.E_i) = \sum_{i = 1}^d (D.E_i)^2$. Taking $N
= \sum_{i = 1}^d (D.E_i)^2$, the claim follows by induction on $a$,
the base case $a = 0$ being trivial.

Now let $D \in N^1(X)$ be ample.  Using the claim of the first part
of the proof we calculate
\begin{gather*}
(D_n.D_n) \ = \  \sum_{i = 0}^{n-1} \sum_{\ell = 0}^{n-1}
(P^iD.P^{\ell}D) \ = \ \sum_{i = 0}^{n-1} (P^iD.P^iD) +  \sum_{0
\leq i < \ell \leq n-1}
2 (P^iD.P^{\ell}D) \\
= \ n(D.D) + \sum_{i=0}^{n-1} iN +  \sum_{m = 1}^{n-1} 2(n-m) (D.P^m
D) + \sum_{m = 1}^{n-1} (n-m-1)(n-m)N.
\end{gather*}

Obviously $n(D.D) \sim n$ and $\sum_{i=0}^{n-1} iN \sim Nn^2$. Using
the result of Lemma~\ref{form-growth-lem}, we get the growth
estimate
\[
\sum_{m = 1}^{n-1} 2(n-m) (D.P^m D) \sim \sum_{m = 1}^{n-1} 2(n-m)
m^j \sim n^{j+2}.
\]
Finally, $\sum_{m = 1}^{n-1} (n-m-1)(n-m)N \sim Nn^3$.  To conclude,
if $\sigma$ is geometric (so $N = 0$) then  $(D_n.D_n) \sim
n^{j+2}$.  If $\sigma$ is non-geometric, then $j = 1$ and $N \neq 0$
and we have  $(D_n.D_n) \sim n^3 = n^{j+2}$ in this case as well.
\end{proof}

We now calculate the growth of the pieces of $\wt{B}(X, \mc{L},
\sigma)$ (for sufficiently positive $\mc{L}$).  We will see later
that this growth is exponential if $\rho > 1$; here we just consider
the case $\rho = 1$.

\begin{proposition}
\label{growth-Ln-prop} Assume Hypothesis~\ref{standard-hyp} and that
$\rho = 1$.  Let $\mc{L} \subseteq \mc{K}$ be a \se\ invertible
sheaf on $X$, and let $\wt{B} = \wt{B}(X, \mc{L}, \sigma)$. Then
$\dim_k \wt{B}_n \sim n^{j+2}$.
\end{proposition}
\begin{proof}
We have $\wt{B}_n \cong  \HB^0(X, \mc{L}_n)$ by definition and
$\HB^i(X, \mc{L}_n) = 0$ for $n  \geq 1$ and $i > 0$ by the \se\
hypothesis. Writing $\mc{L} = \mc{O}_X(D)$, we have $\mc{L}_n =
\mc{O}_X(D_n)$ where $D_n = \sum_{i = 0}^{n-1} P^i(D)$. By the
Riemann-Roch formula, for $n \geq 1$ we have
\[
\dim_k \HB^0(X, \mc{L}_n) = (D_n.D_n)/2 - (D_n.K)/2 + 1 + p_a
\]
where $K$ is the canonical divisor on $X$ and $p_a$ is the
arithmetic genus.   By Lemma~\ref{form-growth-lem2} we have
$(D_n.D_n) \sim n^{j+2}$ since $D$ is ample, and by
Lemma~\ref{form-growth-lem} we have $| (P^n(D).K)|  \preccurlyeq
n^j$ and thus $| (D_n.K) | \preccurlyeq n^{j+1}$.  So $\dim_k
\HB^0(X, \mc{L}_n) \sim n^{j+2}$ as required.
\end{proof}

\section{The lower bound}
\label{lower-sec}

We begin this section by recalling some definitions concerning
graded rings and their growth; see \cite{KL} for more details. Given
any finitely generated $k$-algebra $R$, if $V$ is a
finite-dimensional generating subspace for $R$ containing $1$ then
the \emph{Gelfand-Kirillov dimension} of $R$ is $\GK R =
\overline{\lim} \log_n \dim_k V^n$ (which does not depend on $V$).
The algebra $R$ is said to have \emph{exponential growth} if
$\overline{\lim} (\dim_k V^n)^{1/n}
>1$.  If $R$ is not finitely generated as an algebra, then $\GK R$ is
defined to be the supremum of $\GK R'$ for finitely generated
subalgebras $R' \subseteq R$.

Throughout $\mb{N} = \{0,1,2, \dots \}$ denotes the nonnegative
integers. Let $A$ be an $\mb{N}$-graded $k$-algebra; $A$ is called
\emph{locally finite} if $\dim_k A_n < \infty$ for all $n \geq 0$.
If $A$ is a locally finite $\mb{N}$-graded $k$-algebra which is also
finitely generated as a $k$-algebra, then $\GK A = (\overline{\lim}
\log_n \dim_k A_n) + 1$ (see \cite[Lemma 6.1]{KL}) and moreover $A$
has exponential growth if and only if $\overline{\lim} (\dim_k
A_n)^{1/n} > 1$.

The goal we are heading towards is to understand the growth of
finitely generated $\mb{N}$-graded subalgebras $A \subseteq Q = K[t,
t^{-1}; \sigma]$ where $K/k$ is a finitely generated field extension
of transcendence degree $2$.  An arbitrary such $A \subseteq Q$
could be quite small (for example, if generated as an algebra by a
single element, isomorphic to a polynomial ring in one variable) and
so we concentrate on those $A$ which are large in $Q$ in the sense
of the following definition from \cite{RZ}.
\begin{definition}
\label{bigdef}  A locally finite $\mb{N}$-graded subalgebra $A =
\bigoplus_{n = 0}^{\infty} V_n t^n \subseteq Q = K[t, t^{-1};
\sigma]$ is called a \emph{big} subalgebra of $Q$ if there is some
$n \geq 1$ and an element $u \in V_n$ such that $K$ is the field of
fractions of its subalgebra $k[V_n u^{-1}]$.
\end{definition}

To find a lower bound for the GK-dimension of a big finitely
generated $\mb{N}$-graded algebra $A \subseteq Q = K[t, t^{-1};
\sigma]$, the idea is to compare $A$ with rings of the form
$\wt{B}(X, \mc{L}, \sigma)$, the growth of which we calculated in
the last section. If $\sigma$ is geometric then this is easy and is
already accomplished in \cite[Proposition 5.5(2)]{RZ}: in this case
any big $A$ contains an isomorphic copy of some twisted homogeneous
coordinate ring $B(X, \mc{L}, \sigma)$, and the lower bound is
immediate.   When $\sigma$ is non-geometric, however, there is no
obvious reason that $A$ should contain a copy of some $\wt{B}(X,
\mc{L}, \sigma)$; in fact this seems highly unlikely since $\wt{B}$
is typically infinitely generated. In any case, the same proof as in
\cite[Proposition 5.5(2)]{RZ} does not work and so we use a
different tactic to show that $A$ is not growing more slowly than
$\wt{B}$: we restrict all divisors to a generic curve $E$ on $X$ and
calculate a lower bound for the growth there using some
combinatorial estimates.

First, we need a simple analytic lemma about the growth of functions
satisfying a certain recurrence relation.  Recall the equivalence
relation on growth functions introduced before
Lemma~\ref{rho-facts-lem}.
\begin{lemma}
\label{growth-lem}
 Let $f(n): \mb{N} \to \mb{R}_{> 0}$ be a sequence of positive real numbers which satisfies the
relation $f(n+1) \geq f(n) + f(m(n))$ for some function $m: \mb{N}
\to \mb{N}$ and all $n \gg 0$.
\begin{enumerate}
\item If $m(n) \sim n^{\beta/(\beta+1)}$ for some real $\beta \geq 0$
then $n^{\alpha} \preccurlyeq f(n)$ for every real number $0 <
\alpha < \beta + 1$.
\item If $m(n) \geq n-q$ for some $q \in \mb{N}$ and all $n \gg 0$, then $\delta^n \preccurlyeq
f(n)$ for some real $\delta > 1$.
\end{enumerate}
\end{lemma}
\begin{proof}
(1) Fix some $\alpha$ with $0 < \alpha < \beta +1$.  It is a
calculus exercise to prove that
\[
\lim_{n \to \infty} \frac{(\frac{n+1}{n})^{\alpha} - 1}{\frac{1}{n}}
= \alpha,
\]
and so $((n+1)^{\alpha} - n^{\alpha}) \sim n^{\alpha -1}$.  Thus we
can find a real $\lambda > 0$ such that $(n+1)^{\alpha} \leq
n^{\alpha} + \lambda n^{\alpha -1}$ for all $n \geq n_0$, some $n_0
> 0$. Adjusting $n_0$ higher if necessary, we can find a real
$\epsilon > 0$ such that $m(n) \geq \epsilon n^{\beta/(\beta +1)}$
for all $n \geq n_0$. Again adjusting $n_0$, since $n$ dominates
$n^{\beta/(\beta+1)}$ we can also assume that $n \geq m(n)$ for $n
\geq n_0$.  Finally, the hypothesis $\alpha  < \beta + 1$ implies
that $(\alpha -1) < \alpha \beta/(\beta+1)$ and so with a final
enlargement of $n_0$ we may also assume that $\epsilon^{\alpha}
n^{\alpha \beta/(\beta+1)} \geq \lambda n^{\alpha -1}$ holds for $n
\geq n_0$.

Since $f(n)$ is positive valued, we can choose a real $\gamma > 0$
so that $f(k) \geq \gamma k^{\alpha}$ holds for all $k \leq n_0$.
Let $n \geq n_0$ and assume we have proven already that $f(k) \geq
\gamma k^{\alpha}$ for all $k \leq n$.   Then
\[
f(n+1) \geq f(n) + f(m(n)) \geq \gamma n^{\alpha} + \gamma
m(n)^{\alpha} \geq \gamma n^{\alpha} + \gamma  \epsilon^{\alpha}
n^{\alpha \beta/(\beta+1)}  \geq \gamma n^{\alpha} + \gamma \lambda
n^{\alpha -1} \geq \gamma (n+1)^{\alpha}.
\]
Thus the inequality $f(n) \geq \gamma n^{\alpha}$ holds for all $n$
by induction.

(2) The linear recurrence relation $g(n+1) = g(n) + g(n-q)$ is well
known to give a function $g(n)$ with exponential growth, and the
hypothesis in this case implies that $f(n)$ is growing at least this
fast.
\end{proof}

Next, fix a nonsingular projective curve $C$ over $k$, and let $F =
k(C)$ be its field of rational functions.  A finite dimensional
subspace $W \subseteq F$ generates an invertible sheaf $\mc{M} = W
\mc{O}_C$ on $C$.   We will use the following convenient notation in
the next result: for any $f \in F$ we let $d(f) = \deg W \mc{O}_C$
for $W = k + k f$.  Alternatively, if $(f)$ is the principal divisor
of the rational function $f$ and we write $(f) = P - Q$ where $P$ is
the divisor of zeroes and $Q$ is the divisor of poles, then $d(f) =
\deg P = \deg Q$ since $(k + kf) \mc{O}_C = \mc{O}_C(Q)$. Note that
for an arbitrary vector subspace $W \subseteq F$, $\deg W \mc{O}_C$
gives us no information about $\dim_k W$; in fact, every effective
invertible sheaf on $C$ can be generated by at most $2$ sections. On
the other hand, the next lemma shows that if we have a collection of
subspaces $V_0, V_1, V_2 \dots \subseteq F$, then knowledge of the
growth of $\deg V_i \mc{O}_C$ does allow us to determine the growth
of  $\dim_k V_0 V_1 \dots V_{n-1}$.

\begin{lemma}
\label{curve-growth-lem} Let $C$ be a nonsingular projective curve
over $k$ with $F = k(C)$.  Let $V_0, V_1, V_2, \dots \subseteq F$ be
a sequence of nonzero finite dimensional $k$-vector subspaces of
$F$, and put $d_n = \deg V_n \mc{O}_C$; assume that $d_n
> 0$ for $n \gg 0$.  Let $W_n = V_0 V_1 \dots V_{n-1}$ for each $n$,
and put $e_n = \dim_k W_n$.
\begin{enumerate}
\item Suppose that $d_n \sim n^j$ for some integer $j \geq 0$.  Then for
every real number $0 < \alpha < j+1$, we have $n^{\alpha}
\preccurlyeq e_n \preccurlyeq n^{j+1}$.
\item Suppose that $d_n \sim n^j \rho^n$ for some integer $j \geq 0$ and real $\rho > 1$.  Then $
\delta^n \preccurlyeq e_n$ for some real $\delta > 1$.
\end{enumerate}
\end{lemma}
\begin{proof}
First we give an upper bound on the growth in case (1), so suppose
that $d_n \sim n^j$.  Then $\mc{M}_n = W_n \mc{O}_C$ is an
invertible sheaf of degree $h_n = \sum_{i = 0}^{n-1} d_i$ on $C$.
Since $d_n
> 0$ for $n \gg 0$, $h_{n+1} > h_n$ for $n \gg 0$. In particular,
$h_n > 2g -2$ for $n \gg 0$, where $g$ is the genus of $C$.  Then
$\mc{M}_n$ is nonspecial for $n \gg 0$ (see \cite[Example
IV.1.3.4]{Ha}) and so by the Riemann-Roch formula, $e_n = \dim_k W_n
\leq \dim_k \HB^0(\mc{M}_n) = h_n + 1 -g$.  Thus $e_n \preccurlyeq
h_n \sim n^{j+1}$.

Now we work on the lower bounds.   Replacing each $V_i$ with some
$g_i V_i$ with $0 \neq g_i \in F$ does not affect the numbers $d_i$
or $e_i$, so by making such replacements we may assume that $1 \in
V_i$ for all $i \geq 0$.  Then $W_n \subseteq W_{n+1}$ for all $n$.
Also, since $\dim_k V = 1$ implies that $\deg V \mc{O}_C = 0$, for
$n \gg 0$ we must have $\dim_k V_n \geq 2$. Consider some such $n$.
There is some finite set of points $S \subseteq C$ such that all
elements of $W_n = V_0 V_1 \dots V_{n-1}$ have poles only along
points of $S$. Since $V_n$ is at least $2$-dimensional and contains
$1$, we may choose a generic $f \in V_n$ which does not have zeroes
at any point in $S$, and such that $d(f) = d_n = \deg V_n \mc{O}_C$.

Given any $g \in W_n$, write $(g) = P - Q$ where $P$ is the divisor
of zeroes and $Q$ the divisor of poles of $g$.  Similarly, write
$(f) = P' - Q'$ and $(fg) = (f) + (g) = P'' - Q''$.  By
construction, $P'$ and $Q$ have disjoint support and so $\deg P''
\geq \deg P'$. We conclude that $d(fg) \geq d(f) = d_n$ for all $g
\in W_n$.  Now let $m = m(n)$ be the largest integer such that $\deg
W_m \mc{O}_C < d_n$ (this makes sense, since $\deg W_m \mc{O}_C$ is
eventually strictly monotonic increasing by the first paragraph of
the proof.) Then $d(g') < d_n$ for all $g' \in W_m$. We conclude
that $W_m \cap W_n f = 0$. Since $W_m + W_nf \subseteq W_{n+1}$,
this leads to the formula $e_{n+1} \geq e_n + e_m$.

Now in case (1), we have $d_n \sim n^j$ and $\deg W_n = h_n =
\sum_{i = 0}^{n-1} d_i \sim n^{j+1}$, so it is easy to see that
$m(n) \sim n^{j/(j+1)}$.   Then for every $0 < \alpha < j+1$,
$n^{\alpha} \preccurlyeq e_n$ by Lemma~\ref{growth-lem}(1).

In case (2), suppose that $\mu n^j \rho^n \leq d_n \leq \epsilon n^j
\rho^n$ for all $n \gg 0$, some $0 < \mu < \epsilon$.   Note that
$\deg W_n \mc{O}_C = \sum_{i = 0}^{n-1} d_i \leq \epsilon \sum_{i =
0}^{n-1} i^j \rho^i \leq \epsilon n^j (\rho^n - 1)/(\rho -1)$.  Now
it is easy to see that there is an integer $q > 0$ such that for all
$n \geq 0$, and for any $0 \leq \ell \leq n -q$, we have $\epsilon
\ell^j (\rho^\ell- 1)/(\rho -1) < \mu n^j \rho^n$.  Thus $m(n) \geq
n-q$, and so $\delta^n \preccurlyeq e_n$ for some real $\delta > 1$
by Lemma~\ref{growth-lem}(2).
\end{proof}

Now we are ready to give a lower bound on the growth of certain
subalgebras of $\wt{B}(X, \mc{L}, \sigma)$.
\begin{proposition}
\label{lower-bound-prop}  Assume Notation~\ref{standard-not}, so
$\sigma: X \dra X$ is a stable birational map of a nonsingular
projective model $X$ of $K$, with corresponding growth data $(\rho,
j)$.  Assume that the base field $k$ is uncountable.  By
Lemma~\ref{stab-eff-lem}, pick a very ample invertible sheaf $\mc{L}
\subseteq \mc{K}$ which is \se, and let $U = \HB^0(X, \mc{L})
\subseteq K$. Let
\[
A = k \langle Ut \rangle \subseteq \wt{B}(X, \mc{L}, \sigma) =
\bigoplus_{n \geq 0} \HB^0(X, \mc{L}_n) t^n \subseteq Q = K[t,
t^{-1}; \sigma].
\]
Then $A$ has exponential growth if $\rho >1$, and $\GK A \geq j + 3$
if $\rho = 1$.
\end{proposition}
\begin{proof}
For each $n \in \mb{Z}$ the birational map $\sigma^n: X \dra X$ is
defined except at a finite set of points $S_n$ of $X$.  We may
adjust the embedding $\mc{L} \subseteq \mc{K}$ if necessary so that
$1 \in U = \HB^0(X, \mc{L}) \subseteq K$.  Let $\mc{L} =
\mc{O}_X(D)$, and consider the complete linear system $|D|$ on $X$.
Using Bertini's Theorem and the fact that $k$ is uncountable, we may
choose a Weil divisor $E \in |D|$ with the following properties: (i)
$E$ is an irreducible nonsingular curve, and (ii) $E$ contains none
of the countably many points in $S = \bigcup_{n \geq 1} S_n$.

Let $\mc{I} = \mc{O}_X(-E)$ be the ideal sheaf of the divisor $E$.
It follows from the definition of $\wt{B} = \wt{B}(X, \mc{L},
\sigma)$ that $I = \bigoplus_{n \geq 0} \HB^0(X, \mc{I} \otimes
\mc{L}_n)t^n $ is a homogeneous right ideal of $\wt{B}$. Recalling
the notational convention $U^{\sigma} = \sigma(U)$ from
Section~\ref{growth-sec}, write $U_n = U U^{\sigma} \dots
U^{\sigma^{n-1}} \subseteq K$ for each $n \geq 1$ and $U_0 = k$, so
that $A = \bigoplus_{n \geq 0} U_n t^n$.  For each $i$, the vector
space $U^{\sigma^i}$ generates the sheaf $\mc{L}^{\sigma^i}$ except
possibly at points in $S_i \subseteq S$, by
Lemma~\ref{pullback-sec-lem}. Put $\mc{M}^{(i)} = \mc{L}^{\sigma^i}
\vert_E$ and let $V_i \subseteq \HB^0(E, \mc{M}^{(i)})$ be the image
of $U^{\sigma^i} \subseteq \HB^0(X, \mc{L}^{\sigma^i})$ under the
restriction of sections map. By the choice of $E$, $\mc{M}^{(i)}$ is
an invertible sheaf on $E$ which is generated by the sections in
$V_i$ everywhere. Similarly, put $\mc{M}_n = \mc{L}_n \vert_E$ and
let $W_n \subseteq \HB^0(E, \mc{M}_n)$ be the image of $U_n
\subseteq \HB^0(X, \mc{L}_n)$ under the restriction of sections map;
again, $W_n$ generates $\mc{M}_n$ on $E$. Choose arbitrary
embeddings $\mc{M}^{(i)} \subseteq \mc{F}$ for each $i$, where
$\mc{F}$ is the constant sheaf on $F = k(E)$, which determines
embeddings of $\mc{M}_n = \mc{M}^{(0)} \otimes \mc{M}^{(1)} \otimes
\dots \otimes \mc{M}^{(n-1)} \subseteq \mc{F}$ for each $n$. This
also embeds $V_i \subseteq F$ for each $i$ and $W_n \subseteq F$ for
each $n$, and since $U_n = U U^{\sigma} \cdots U^{\sigma^{n-1}}$ in
$K$, we have $W_n = V_0 V_1 \cdots V_{n-1}$ in $F$. Observe that
there is a graded vector space map $\theta: A \to \bigoplus_{n \geq
0} \HB^0(E, \mc{M}_n)$ given by restriction of sections, with kernel
$A \cap I = \bigoplus_{n \geq 0} (U_n \cap \HB^0(X, \mc{I} \otimes
\mc{L}_n))t^n $ and image $\bigoplus_{n \geq 0} W_n$.  (In fact,
$\bigoplus_{n \geq 0} \HB^0(E, \mc{M}_n)$ is naturally a graded
right $A$-module and $\theta$ is a right module map, but we won't
need this.)

In case  $\rho > 1$, to show that $A$ has exponential growth it will
clearly suffice to show that $\dim_k W_n$ grows exponentially with
$n$. In case $\rho = 1$, note that by construction $A_1 = \HB^0(X,
\mc{L})t$ and $\mc{L}$ contains a nonzero section vanishing along
the curve $E$, so $(A \cap I)_1 \neq 0$. In particular, the
GK-dimension of the right $A$-module $A/(A \cap I)$ satisfies $\GK
A/(A \cap I) \leq \GK A -1$ by \cite[Proposition 5.1(e)]{KL}.  Then
if we can show that $n^{\alpha} \preccurlyeq \dim_k W_n$ holds for
all $\alpha < j +1$, we will have $j + 2 \leq \GK A/(A \cap I)$ by
\cite[Lemma 6.1(b)]{KL}, and so $j + 3 \leq \GK A$.

Now on the curve $E$, putting $d_n = \deg V_n \mc{O}_E = \deg
\mc{M}^{(n)}$ we have $d_n = ((\sigma^n)^*(E).E)$ by \cite[Lemma
V.1.3]{Ha}; then since $E$ is ample, $d_n \sim n^j \rho^n$ by
Lemma~\ref{form-growth-lem}.  In particular, notice that $d_n > 0$
for all $n \geq 0$.   Finally, by Lemma~\ref{curve-growth-lem},
putting $e_n = \dim_k W_n$ then $e_n$ has exponential growth if
$\rho
 > 1$, while $n^{\alpha} \preccurlyeq e_n$ for all $\alpha < j + 1$ if
$\rho = 1$.  The result follows.
\end{proof}

\section{Proof of the main theorem}
\label{final-sec}

We now put together the various estimates already proved to
calculate the GK-dimension of big subalgebras $A \subseteq K[t,
t^{-1}; \sigma]$.
\begin{theorem}
\label{final-thm} Let $k$ be an uncountable algebraically closed
field. Let $Q = K[t, t^{-1}; \sigma]$, where $K$ is a finitely
generated field extension of $k$ with $\trdeg K/k = 2$ and $\sigma
\in \aut_k(K)$.  Let $(\rho, j)$ be the growth data associated to
$\sigma$, as in Notation~\ref{standard-not}.   Given any big locally
finite $\mb{N}$-graded subalgebra $A \subseteq Q$, then $\GK A$ is
determined as follows:
\begin{enumerate}
\item If $\rho > 1$, then $A$ has exponential growth (in particular, $\GK A =
\infty$).
\item If $\rho = 1$, then $\GK A = j  + 3$.  In particular, $\GK A
\in \{3,4,5\}$ and $\GK A = 4$ if and only if $\sigma$ is
non-geometric.
\end{enumerate}
\end{theorem}
\begin{proof}
Suppose first that $\rho =1$.  Choose a nonsingular projective model
$X$ of $K$ so that the birational map $\sigma: X \dra X$ satisfies
Hypothesis~\ref{standard-hyp}.  Write $A = \bigoplus_{n \geq 0} U_n
t^n$ where $U_n \subseteq K$. Assume for the moment that $A$ is
finitely generated as an algebra. Since $A_0$ is a
finite-dimensional $k$-algebra which is a domain and $k$ is
algebraically closed, $A_0 = k$. Suppose that $A$ is generated by
elements of degree $\leq d$. Letting $W = k + U_1 + \dots U_d$, then
$A$ is contained in the ring $k \langle W t \rangle \subseteq Q$.
Moreover, we can find a very ample invertible sheaf $\mc{L}$ with an
embedding $\mc{L} \subseteq \mc{K}$ in the constant sheaf of
rational functions such that $W \subseteq \HB^0(X, \mc{L}) \subseteq
K$ \cite[Lemma 5.2]{RZ}. Replace $\mc{L}$ by a power
$\mc{L}^{\otimes m}$ if necessary, so that $\mc{L}$ is \se, using
Lemma~\ref{stab-eff-lem}.  Then $k \langle Wt \rangle \subseteq
\bigoplus_{n \geq 0} \HB^0(X, \mc{L}_n)t^n = \wt{B}(X, \mc{L},
\sigma) \subseteq Q$.  But in Proposition~\ref{growth-Ln-prop} we
saw that $\dim_k \wt{B}_n \sim n^{j+2}$; hence $\GK A \leq j + 3$.
Now since this estimate holds for all finitely generated $A$, $\GK A
\leq j + 3$ holds for an arbitrary locally finite $\mb{N}$-graded
subalgebra $A$ of $Q$, by definition.

Now we find a lower bound for $\GK A$.  Choose any very ample
invertible sheaf $\mc{L} \subseteq \mc{K}$ which is \se\ and let $W
= \HB^0(X, \mc{L}) \subseteq K$. Suppose first that $\rho
> 1$, so we want to prove that $A$ has exponential growth. In this
case, by \cite[Proposition 1.4]{RZ} it suffices to find any finitely
generated $\mb{N}$-graded subalgebra of $Q$ with exponential growth,
and so we will prove that $k \langle Wt \rangle \subseteq Q$ has
exponential growth.

If instead $\rho = 1$, then the first part of the proof shows that
$\GK A < \infty$, so $A$ is an Ore domain by \cite[Proposition
4.13]{KL}. Then again setting $A = \bigoplus_{n \geq 0} U_n t^n$,
the same argument as in \cite[Lemma 5.3]{RZ} (using that $A$ is big
in $Q$) shows that there is some $n > 0$ and $z \in U_n$ such that
$Wz \subseteq U_n$. Setting $t' = t^n$ and $\sigma' = \sigma^n$, the
Veronese ring $A' = A^{(n)}$ contains $k \langle Wz t' \rangle$,
which is a big subalgebra of $Q' = Q^{(n)} = K[t', (t')^{-1};
\sigma']$. Note that the induced map $\sigma' = \sigma^n: X \dra X$
is again a stable birational map, with the same associated growth
data $(\rho, j) = (1,j)$. Moreover, $k \langle Wz t' \rangle \cong k
\langle W t' \rangle$. It is enough to put a lower bound on $\GK
A'$, so we now change notation back by removing the primes, and our
task is again to find a lower bound for the GK-dimension of $k
\langle W t \rangle \subseteq Q$.

Thus in both cases for $\rho$, the rings $k \langle Wt \rangle
\subseteq \bigoplus_{n \geq 0} \HB^0(X, \mc{L}_n)t^n  = \wt{B}(X,
\mc{L}, \sigma)$ now satisfy the hypothesis of
Proposition~\ref{lower-bound-prop}. By that proposition, $k \langle
Wt \rangle$ has exponential growth if $\rho
> 1$ and $\GK k \langle Wt \rangle \geq j + 3$ if $\rho = 1$.

Taking the upper and lower bounds together, the calculation of $\GK
A$ is complete.  The rest of part (2)  follows immediately from the
classification in Theorem~\ref{DF-sum-thm}.
\end{proof}

\emph{Proof of Theorem~\ref{intro-thm}.}  The theorem follows
immediately from Theorem~\ref{final-thm} if $k$ is uncountable, so
we just need to reduce to this case.   Suppose that $k$ is any
algebraically closed field.  Let $A \subseteq K[t, t^{-1}; \sigma]$
be a big locally finite $\mb{N}$-graded subalgebra, where $K/k$ is a
finitely generated field extension with $\trdeg K/k = 2$.  Choose
any field extension $k \subseteq L$ where $L$ is uncountable and
algebraically closed. As was also noted in
Lemma~\ref{base-change-lem}, since $k$ is algebraically closed, $K
\otimes_k L$ is again a commutative domain; letting $F$ be its field
of fractions, $\sigma$ extends uniquely to an automorphism
$\wt{\sigma} \in \aut_L F$. Consider $\wt{A} = A \otimes_k L
\subseteq K [t, t^{-1}; \sigma] \otimes_k L \subseteq F[t, t^{-1};
\wt{\sigma}]$.  It is easy to see that $\wt{A}$ is a big finitely
$\mb{N}$-graded subalgebra of $F[t, t^{-1}; \wt{\sigma}]$, where
again $F/L$ is a finitely generated field extension with $\trdeg F/L
= 2$.  Also, $\GK_L \wt{A} = \GK_k A$. Now the result follows from
Theorem~\ref{final-thm}, together with Lemma~\ref{base-change-lem}.
\hfill $\Box$

\section*{Acknowledgments}
This paper benefited from conversations with James Zhang, Mike
Artin, and Toby Stafford, and we thank them.  We would also like to
thank the referee for a very careful reading of the manuscript.

\providecommand{\bysame}{\leavevmode\hbox
to3em{\hrulefill}\thinspace}
\providecommand{\MR}{\relax\ifhmode\unskip\space\fi MR }
\providecommand{\MRhref}[2]{%
  \href{http://www.ams.org/mathscinet-getitem?mr=#1}{#2}
} \providecommand{\href}[2]{#2}

\end{document}